\documentclass[noinfoline,11pt]{imsart}

\RequirePackage[OT1]{fontenc}
\RequirePackage{amsthm,amsmath}
\RequirePackage[numbers]{natbib}
\RequirePackage[colorlinks,citecolor=blue,urlcolor=blue]{hyperref}
\usepackage{graphicx}

\startlocaldefs
\usepackage{amsfonts,amstext, amssymb,amsthm, amsmath, rotating, latexsym,wasysym}
\usepackage[margin=1in]{geometry}
\usepackage{caption}
\usepackage{enumerate}
\usepackage{hyperref}
\usepackage{appendix}
\usepackage{booktabs}
\usepackage{graphicx}
\usepackage{caption}
\usepackage{float}
\usepackage[table,xcdraw]{xcolor}
\usepackage{relsize}
\usepackage[linesnumbered,ruled,vlined]{algorithm2e}
\usepackage{color}
\usepackage{mathtools}
\mathtoolsset{showonlyrefs}
\usepackage{graphicx,bm,bbm,bigints}
\usepackage{multirow}
\usepackage{booktabs}
\usepackage{mathrsfs}
\usepackage{enumitem}
\usepackage[caption=false]{subfig}
\usepackage[countmax]{subfloat}
\usepackage{chngcntr}
\usepackage{apptools}
\usepackage{footnote}
\numberwithin{equation}{section}
\theoremstyle{plain}

\newcommand{\y}{\mathbf{y}}
\newcommand{\bb}{\boldsymbol{\beta}}

\newcommand{\bl}{\boldsymbol{\lambda}}

\newtheorem{theorem}{Theorem}
\newtheorem{lemma}{Lemma}
\newtheorem{corollary}{Corollary}

\AtAppendix{\counterwithin{lemma}{section}}
\AtAppendix{\counterwithin{result}{section}}
\AtAppendix{\counterwithin{corollary}{section}}
\AtAppendix{\counterwithin{theorem}{section}}

\makeatletter

\newcommand{\Rome}[1]{\expandafter\@slowromancap\romannumeral #1@}
\makeatother
\newcommand{\be}{\begin{equation}}
\newcommand{\ee}{\end{equation}}
\newcommand{\ba}{\begin{eqnarray}}
\newcommand{\ea}{\end{eqnarray}}
\newcommand{\bee}{\begin{equation*}}
\newcommand{\eee}{\end{equation*}}
\newcommand{\baa}{\begin{eqnarray*}}
	\newcommand{\eaa}{\end{eqnarray*}}

\newcommand{\R}{\mathbb{R}}

\definecolor{lgray}{gray}{0.70}
\newcommand{\rpm}{\sbox0{$1$}\sbox2{$\scriptstyle\pm$}
	\raise\dimexpr(\ht0-\ht2)/2\relax\box2 }

\newcount\colveccount
\newcommand*\colvec[1]{
	\global\colveccount#1
	\begin{pmatrix}
		\colvecnext
	}
	\def\colvecnext#1{
		#1
		\global\advance\colveccount-1
		\ifnum\colveccount>0
		\\
		\expandafter\colvecnext
		\else
	\end{pmatrix}
	\fi
}

\makeatletter
\newcommand\@erelb@r[1]{%
	\mathrel{\tikz[baseline=-.5ex]\draw[#1] (0,0)--(0.3,0);}
}
\newcommand{\erelbar}[1]{\@erelbar#1}
\def\@erelbar#1#2{%
	\ifcase\numexpr#1*4+#2\relax
	\@erelb@r{-}\or % 00
	\@erelb@r{->}\or % 01
	\@erelb@r{-|}\or % 02
	\@erelb@r{->|}\or % 03
	\@erelb@r{<-}\or % 10
	\@erelb@r{<->}\or % 11
	\@erelb@r{<-|}\or % 12
	\@erelb@r{<->}\or % 13
	\@erelb@r{|-}\or % 20
	\@erelb@r{|->}\or % 21
	\@erelb@r{|-|}\or % 22
	\@erelb@r{|<->|}\or % 23
	\@erelb@r{|<-}\or % 30
	\@erelb@r{|<->}\or % 31
	\@erelb@r{|<-|}\or % 32
	\@erelb@r{|<->|} % 33
	\else
	\@wrong
	\fi
}
\makeatother

\makesavenoteenv{tabular}

\allowdisplaybreaks

\endlocaldefs

\begin{document}
  
\begin{frontmatter}

\title{Convergence properties of data augmentation algorithms for high-dimensional robit regression}
\runtitle{Convergence of robit regression DA algorithm}

\begin{aug}
\author[a]{\fnms{Sourav} \snm{Mukherjee} \ead[label=e1]{souravmukherjee@ufl.edu}},
\author[a]{\fnms{Kshitij} \snm{Khare}\ead[label=e2]{kdkhare@stat.ufl.edu}}
\and 
\author[b]{\fnms{Saptarshi} \snm{Chakraborty}\ead[label=e3]{chakrab2@buffalo.edu}}
\address[a]{Department of Statistics, University of Florida}
\address[b]{Department of Biostatistics, State University of New York at 
Buffalo}

\printead{e1,e2,e3}\\

\runauthor{S. Mukherjee, K. Khare, and S. Chakraborty}
\end{aug}

\begin{abstract}
    The logistic and probit link functions are the most common choices for regression models with a binary response. However, these choices are not robust to the presence of outliers/unexpected observations. The robit link function, which is equal to the inverse CDF of the Student's $t$-distribution, provides a robust alternative to the probit and logistic link functions. A multivariate normal prior for the regression coefficients is the standard choice for Bayesian inference in robit regression models. The resulting posterior density is  intractable and a Data Augmentation (DA) Markov chain is used to generate approximate samples from the desired posterior distribution. Establishing geometric ergodicity for this DA Markov chain is important as it provides theoretical guarantees for asymptotic validity of MCMC standard errors for desired posterior expectations/quantiles. Previous work \cite{Roy:2012:robit} established geometric ergodicity of this robit DA Markov chain assuming (i) the sample size $n$ dominates the number of predictors $p$, and (ii) an additional constraint which requires the sample size to be bounded above by a fixed constant which depends on the design matrix $X$. In particular, modern high-dimensional settings where $n < p$ are not considered. In this work, we show that the robit DA Markov chain is trace-class (i.e., the eigenvalues of the corresponding Markov operator are summable) for arbitrary choices of  the sample size $n$, the number of predictors $p$, the design matrix $X$, and the prior mean and variance parameters. The trace-class property implies geometric ergodicity. Moreover, this property allows us to conclude that the sandwich robit chain (obtained by inserting an inexpensive extra step in between the two steps of the DA chain) is strictly better than the robit DA chain in an appropriate sense, and enables the use of recent methods to estimate the spectral gap of trace class DA Markov chains.  
\end{abstract}

\begin{keyword}[class=MSC]
\kwd[Primary ]{60J05}
\kwd{60J20}
\end{keyword}

\begin{keyword}
\kwd{Markov chain Monte Carlo}
\kwd{geometric ergodicity}
\kwd{High-dimensional binary regression}
\kwd{robust model}
\kwd{trace class}
\end{keyword}

\tableofcontents

\end{frontmatter}

\section{Introduction}

\noindent
Consider a regression setting with $n$ independent binary responses $Y_1, 
Y_2, \cdots, Y_n$ and corresponding predictor vectors ${\bf x}_1, 
{\bf x}_2, \cdots, {\bf x}_n \in \mathbb{R}^p$, such that 
$$
P(Y_i = 1 \mid \bb) = F({\bf x}_i^T \bb) 
$$

\noindent
for $1 \leq i \leq n$. Here $F$ is a strictly increasing cumulative 
distribution function, and $G = F^{-1}$ is referred to as the link 
function. Two popular choices of $F$ are given by $F(x) = 
\frac{e^x}{1+e^x}$ (logistic link) and $F(x) = \Phi(x)$ (probit link) where
$\Phi(x)$ denotes the standard normal CDF. However, it is well known that 
estimates of $\bb$ produced for both these choices are not robust to 
outliers \cite{Pregibon:1982, GH:2007}. To address such settings, $F$ is 
set to be the CDF of the Student's $t$-distribution, and the corresponding 
model is referred to as the robit regression model 
\cite{Chuanhai:2004}. For binary responses, 
an outlier is an unexpected observation with large value(s) of the 
predictor(s) and a misclassified response, and the robit model effectively 
down-weights such outliers to produce a better fit \cite{GHV:2020}. 

Following \cite{Roy:2012:robit, Albert:Chib:1993} we consider a Bayesian 
robit regression model specified as follows. 
\begin{eqnarray}
& & P(Y_i = 1 \mid \bb) = F_\nu ({\bf x}_i^T \bb) \mbox{ for } 1 \leq i 
\leq n \nonumber\\
& & \bb \sim \mathcal{N}_p \left( \bb_a, \Sigma_a^{-1} \right), 
\label{robitmodel}
\end{eqnarray}

\noindent
where $F_{\nu}$ denotes the CDF of the Student's $t$-distribution with 
$\nu$ degrees of freedom (with location $0$, scale $1$) and $\mathcal{N}_p$
denotes the $p$ variate normal distribution. Let ${\bf Y}$ denote the 
response vector, and let $\pi(\bb \mid {\bf y})$ denote the posterior 
density of $\bb$ given ${\bf Y} = {\bf y}$. As demonstrated in 
\cite{Roy:2012:robit}, the posterior density $\pi(\bb \mid {\bf y})$ is 
intractable in the sense that relevant posterior expectations are ratios of
two intractable integrals and are not available in closed form. Also, 
generating IID samples from this density is computationally feasible even 
for moderately large values of $p$. To resolve this, 
\cite{Roy:2012:robit} develops a clever and effective Data Augmentation 
(DA) approach which can be used to construct a computationally tractable 
Markov chain which has $\pi(\bb \mid {\bf y})$ as its stationary density. 
We describe this Markov chain below. 

Let $t_\nu (\mu, \sigma)$ denote the Student's $t$-distribution with $\nu$ 
degrees of freedom, location $\mu$ and scale $\sigma$. Consider unobserved 
latent variables $(Z_1, \lambda_1), (Z_2, \lambda_2), \cdots, (Z_n, 
\lambda_n)$ which are mutually independent and satisfy $Z_i \mid \lambda_i 
\sim \mathcal{N}({\bf x}_i^T \bb, 1/{\bf \lambda}_i)$ and $\lambda_i \sim 
\mbox{Gamma}(\nu/2,\nu/2)$. Then, it can be shown that the marginal 
distribution of $Z_i$ is given by $Z_i \sim t_\nu ({\bf x}_i^T\bb, 1)$. If 
$Y_i$ is defined as the indicator of $Z_i$ taking positive values, i.e., 
$Y_i = 1_{\{Z_i > 0\}}$, then $P(Y_i = 1 \mid \bb) = P(z_i > 0) = F_\nu 
({\bf x}_i^T \bb)$, which is consistent with the robit regression model 
specified in (\ref{robitmodel}). Straightforward calculations (see 
\cite{Roy:2012:robit}) now show the following. 
\begin{itemize}
    \item $(Z_1, \lambda_1), (Z_2, \lambda_2), \cdots, (Z_n, 
    \lambda_n)$ are conditionally independent given $\bb, {\bf Y} = {\bf y}$. Also, 
    \begin{eqnarray}
    & & Z_i \mid \bb, {\bf y} \sim Tt_\nu ({\bf x}_i^T \bb, y_i) \label{zdist}\\
    & & \lambda_i \mid {Z_i = z_i}, \bb, {\bf y} \sim \mbox{Gamma} \left( 
    \frac{\nu+1}{2}, \frac{\nu + (z_i - {\bf x}_i^T \bb)^2}{2} \right). 
    \label{lambdadist}
    \end{eqnarray}
    
    \noindent
    Here $Tt_v ({\bf x}_i^T \bb, y_i)$ denotes the $t_\nu$ distribution with location ${\bf x}_i^T \bb$ and scale $1$, truncated to $\mathbb{R}_+$ if $y_i = 1$ and to $\mathbb{R}_-$ if $y_i = 0$. 
    
    \item Let ${\bf Z} = (Z_1, \cdots, Z_n)^T$, $\bl = (\lambda_1, 
    \cdots, \lambda_n)^T$ and $\Lambda$ denote a diagonal matrix whose diagonal is given by entries of $\bl$. The conditional distribution of $\bb$ given ${\bf Z} = {\bf z}, \bl, {\bf Y} = \y$ is $\mathcal{N}_p \left( (X^T \Lambda X + \Sigma_a)^{-1} (X^T \Lambda {\bf z} + \Sigma_a \bb_a) , \ (X^T \Lambda X + \Sigma_a)^{-1} \right)$. 
\end{itemize}

These observations are used in \cite{Roy:2012:robit} to construct a DA 
Markov chain $\{\bb^{(m)}\}_{m \geq 0}$ on $\mathbb{R}^p$ (with stationary 
density $\pi (\bb \mid \y)$) whose one step transition from $\bb^{(m)}$ to 
$\bb^{(m+1)}$ is  described in the following Algorithm~\ref{algo:robit_da}.
\begin{algorithm}
\caption{$(m+1)$-st Iteration of the Robit Data Augmentation Algorithm}
\label{algo:robit_da}
\begin{enumerate}
    \item Make independent draws from the $Tt_\nu ({\bf x}_i^T \bb^{(m)}, y_i)$ distributions for $1 \leq i \leq n$. Denote the respective draws by $z_1, z_2, \cdots, z_n$. Draw $\lambda_i$ from the Gamma$\left(\frac{\nu+1}{2}, \frac{\nu + \left(z_i - {\bf x}_i^T \bb^{(m)}\right)^2}{2} \right)$ distribution. 
    \item Draw $\bb^{(m+1)}$ from the $\mathcal{N}_p \left( (X^T \Lambda X + \Sigma_a)^{-1} (X^T \Lambda {\bf z} + \Sigma_a \bb_a), \ (X^T \Lambda X + \Sigma_a)^{-1} \right)$ distribution.
\end{enumerate}
\end{algorithm}
\noindent
Harris ergodicity of the robit DA Markov chain $\{\bb^{(m)}\}_{m \geq 0}$ 
obtained through Algorithm~\ref{algo:robit_da} is established in \cite{Roy:2012:robit}. Suppose 
a posterior expectation $E_{\pi(\cdot \mid \y)}[h(\bb)]$ (assumed to exist)
is of interest. Then by Harris ergodicity, the cumulative averages 
$\frac{1}{m+1} \sum_{r=0}^m h(\bb^{(r)})$ converge to $E_{\pi(\cdot \mid 
\y)}[h(\bb)]$ as $m \rightarrow \infty$, and can be used to approximate 
the desired posterior expectation. However, any such approximation is 
useful only with an estimate of the associated error. The standard approach
for obtaining such error estimates is to establish a Markov chain central 
limit theorem (CLT) which guarantees that 
$$
\sqrt{m} \left( \frac{1}{m+1} \sum_{r=0}^m h(\bb^{(r)}) - E_{\pi(\cdot \mid
\y)}[h(\bb)] \right) \stackrel{D}{\rightarrow} \mathcal{N} (0, \sigma_h^2)
$$

\noindent
as $m \rightarrow \infty$, and then construct a consistent estimate 
$\hat{\sigma}_h$ of the asymptotic standard deviation $\sigma_h$. A key 
sufficient condition for establishing a Markov chain CLT is {\it geometric
ergodicity}. A Markov chain is geometrically ergodic if the total variation
distance between its distribution after $m$ steps and the stationary 
distribution converges to $0$ as $m \rightarrow \infty$. To summarize, establishing 
geometric ergodicity of the robit DA chain is critical for obtaining 
asymptotically valid standard errors for Markov chain based estimates of 
posterior quantities. 

With this motivation \cite{Roy:2012:robit} investigated and established geometric 
ergodicity of the robit DA chain. However, it is assumed that the design matrix is full rank (which 
implies $n \geq p$ and rules out high-dimensional settings), $\Sigma_a = g^{-1} X^T X$ and that 
$$
n \leq \frac{g^{-1} \nu}{(\nu+1)(1+2\sqrt{\bb_a^T X^T X \bb_a})}. 
$$

\noindent
The last upper bound on $n$ involving the design matrix, the prior mean and covariance and the 
degrees of freedom $\nu$ is in particular very restrictive. While we found that these conditions can be relaxed to some extent by a tighter drift and 
minorization analysis, the resulting constraints still remain quite restrictive. Geometric ergodicity results for the 
related probit DA chain (see \cite{Chakraborty:Khare:2017}) do not require such assumptions, and quoting 
from \cite[Page 2469]{Roy:2012:robit} ``Ideally, we would like to be able to say that the DA algorithm 
is geometrically ergodic for any $n, \nu, \y, X, \bb_a, \Sigma_a$". 

In the probit DA setting (see \cite{Albert:Chib:1993}) the latent variables $Z_i$ have a normal 
distribution, and there is no need to introduce the additional latent variables $\lambda_i$. This 
additional layer of latent variables creates additional complexity in the structure of the robit 
DA chain which makes the convergence analysis significantly more challenging compared to the 
probit DA chain analyses undertaken in \cite{Roy:Hobert:2007, Chakraborty:Khare:2017}. 

It is not clear if the restrictive conditions listed above for geometric 
ergodicity of the robit DA chain are really necessary or if they are an 
artifact of the standard drift and minorization technique used in 
\cite{Roy:2012:robit} for establishing geometric ergodicity. We take a 
completely different approach, and focus on investigating the {\it trace 
class} property for the robit DA chain. A Markov chain with stationary 
density $\pi$ is trace class if the corresponding Markov operator on $L^2 
(\pi)$ has a countable spectrum and the corresponding eigenvalues are 
summable. The trace class property implies geometric ergodicity, and can be
established by showing that an appropriate integral involving the 
transition density of the Markov chain is finite (see 
Section \ref{sec:trace.class.property}). As the main technical contribution
of this paper, {\it we establish that the robit DA chain is trace class for
any $n, \nu > 2, \y, X, \bb_a, \Sigma_a$}. This in particular {\it establishes 
geometric ergodicity for any $n, \nu > 2, \y, X, \bb_a, \Sigma_a$, and 
significantly generalizes existing results} in \cite{Roy:2012:robit}. 

The trace class property is much stronger than geometric ergodicity, and 
the bounding of the relevant integral can get quite involved and 
challenging (see for example the analysis is \cite{Chakraborty:Khare:2017} 
and in Section \ref{sec:trace.class.property}). It is therefore not 
surprising that the conditions needed for establishing geometric 
ergodicity through the trace class approach have typically been stronger 
than those needed to establish geometric ergodicity using the drift and 
minorization approach (for chains where both such analyses have been 
successful). For example, \cite{Chakraborty:Khare:2017} considers 
convergence analysis of the probit DA chain with a proper prior for $\bb$ 
as in (\ref{robitmodel}). Using drift and minorization geometric ergodicity
is established for all $n, \y, X, \bb_a, \Sigma_a$, but the trace class 
property was only established under some constraints on $X$ and $\Sigma_a$ 
(see \cite[Theorem 2]{Chakraborty:Khare:2017}). Hence, it is quite 
interesting that for the robit DA chain the reverse phenomenon holds: the 
drift and minorization approach needs stronger conditions to succeed than 
the trace class approach. Essentially, the additional layer of latent 
variables $\lambda_i$ introduced in the robit setting severely hampers the 
drift and minorization analysis, but (with careful additional analysis) 
eases the path for showing the finiteness of the relevant trace class 
integral. 

Establising the trace class property of the DA chain gives additional 
benefits on top of geometric ergodicity. Using results from \cite{Khare:Hobert:2011}, 
it can now be concluded that the sandwich robit DA chain constructed in 
\cite{Roy:2012:robit} is also trace class and is strictly better than 
the robit DA chain (in the sense that the spectrum of the latter strictly 
dominates the spectrum of the former). Also, the trace class property is 
a key sufficient condition for using recent approaches in \cite{Chakraborty:Khare:2019, QHK:2019} to 
estimate the spectral gap of Markov chains. The remainder of the paper 
is organized as follows. Section \ref{sec:trace.class.property} contains 
the proof of the trace class property for the robit DA chain. 
Section \ref{sec:numerical:illustration} provides numerical illustrations of various chains on two 
real datasets, one with $n \geq p$, and one with $n < p$. Additional mathematical results needed for the 
proof of the trace class property are provided in an 
appendix. 

%A supplementary document shows that the conditions in 
%\cite{Roy:2012:robit} for geometric ergodicity using 
%the drift and minorization approach can be relaxed. 

\section{Trace-class property for the DA chain}\label{sec:trace.class.property}
Recall from \cite{Roy:2012:robit} that the DA Markov chain has associated transition density given by
\begin{align}
    k\left(\bb, \bb^{'}\right) \ &= \ \int_{\R^n_{+}}\int_{\R^n}\pi\left(\bb^{'} \Big| \bl, {\bf z}, \y \right)\pi\left(\Lambda, {\bf z} \Big| \bb, {\bf y}\right)d {\bf z} \ d\bl\\
    &= \ \int_{\R^n_{+}}\int_{\R^n}\pi \left(\bb^{'} \Big| \bl, {\bf z}, {\bf y}\right)\pi\left(\Lambda \Big| {\bf z}, \bb, {\bf y}\right)\pi\left({\bf z} \Big| \bb, {\bf y}\right)d{\bf z} \ d\bl
    \label{eq:markov.transition.density}
\end{align}
Let $L_0^2(\pi(. | \y))$ denote the space of square-integrable functions with mean zero (with respect to the posterior density $\pi(\bb|\y)$). Let $K$ denote the Markov operator on $L_0^2(\pi(. | \y))$ associated with the transition density $k$. Note that the Markov transition density $k$ is reversible with respect to its invariant distribution, and $K$ is a positive, self-adjoint operator. The operator $K$ is \textit{trace class} (see J\"{o}rgens \cite{jorgens.linear.integral.operator}) if
\begin{equation}
    I := \int_{\R^p}k\left(\bb, \bb\right)d\bb \ < \ \infty
    \label{eq:trace.class.condition}
\end{equation}
If the trace-class property holds, then the spectrum of
$K$ is countable and the corresponding eigenvalues are
summable. This in particular implies that $K$ is 
compact, and the associated Markov chain is 
geometrically ergodic.

The following theorem shows that the Markov 
operator $K$ corresponding to the robit DA 
chain is trace class under very general 
conditions. 

\begin{theorem} \label{thm1}
    For $\nu >2$, the Markov operator $K$ corresponding to the 
    DA Markov chain is trace-class for an arbitrary choice of 
    the design matrix $X$, sample size $n$, number of 
    predictors $p$, prior mean vector $\bb_a$, and (positive 
    definite) prior precision matrix $\Sigma_a$. 
\end{theorem}

\begin{proof}
We shall show that \eqref{eq:trace.class.condition} 
holds for the DA Markov chain. The proof is quite 
lengthy and involved, and we have tried to make it 
accessible to the reader by highlighting the major 
steps/milestones. 

We begin by fixing our notations. Let $I_A(.)$  be the indicator function of the set $A$ and $\phi\left(x; a, b\right)$  be the univariate normal density evaluated at point $x$ with mean $a$ and variance $b$. Further, let $\phi_p\left(\bf{x}; \bm{\mu}, \bm{\Sigma}\right)$ denote the $p$-variate normal density with mean vector ${\bm \mu}$, covariance matrix $\bm{\Sigma}$, evaluated at a vector ${\bf x}\in \R^p$. Finally, let $q\left(\omega; a, b\right) = b^a \omega^{a-1}e^{-b\omega}/\Gamma(a)$ be the gamma density evaluated at  $\omega$ with shape parameter $a$ and rate parameter $b$.

Note from Section $2.1$ of \cite{Roy:2012:robit} that the joint posterior density of $\left(\bb, \bl, {\bf 
z}\right)$ is given by
\begin{align}
    &\pi\left(\bb, \bl, {\bf z} \Big| \y \right)\\
    &= \frac{1}{m(\y)}\left[\prod_{i=1}^{n} \Big\{I_{\R_+}\left(z_i\right)I_{\{1\}}\left(y_i\right) + I_{\R_-}\left(z_i\right)I_{\{0\}}\left(y_i\right)\Big\}\phi\left(z_i; {\bf x}_i^T\bb, \frac{1}{\lambda_i}\right)q\left(\lambda_i; \frac{\nu}{2}, \frac{\nu}{2}\right)\right]\\
    &\qquad\qquad\qquad \times \phi_p\left(\bb; \bb_a, \Sigma_a^{-1}\right)\\
    &= \frac{1}{m(\y)}\Bigg[\prod_{i=1}^{n} \Big\{I_{\R_+}\left(z_i\right)I_{\{1\}}\left(y_i\right) + I_{\R_-}\left(z_i\right)I_{\{0\}}\left(y_i\right)\Big\}\times \frac{\sqrt{\lambda_i}}{\sqrt{2\pi}}\exp\left\{-\frac{\lambda_i}{2}\left(z_i - {\bf x}_i^T\bb\right)^2\right\}\\
    &\qquad\qquad \times \frac{\left(\frac{\nu}{2}\right)^{\frac{\nu}{2}}}{\Gamma\left(\frac{\nu}{2}\right)} \ \lambda_i^{\frac{\nu}{2} - 1}\exp\left\{-\frac{\nu}{2}\lambda_i\right\}\Bigg]\\
    &\qquad\qquad \times \left(2\pi\right)^{-\frac{p}{2}} \sqrt{\text{det}\left(\Sigma_a\right)} \ \exp\left[-\frac{1}{2}\left(\bb - \bb_a\right)^T\Sigma_a\left(\bb - \bb_a\right)\right]
    \label{eq:joint.density.Beta_Lambda_Z.TraceClass}
\end{align}

\noindent
for $\bl \in \R_+^n, \ {\bf z} \in \R^n, \ \bb 
\in \R^p$. 

\medskip

\noindent
{\bf Step I: A useful linear reparametrization to adjust for the prior mean $\bb_a$ and 
derivation of associated conditionals.} 
Consider the following reparametrization
$$\Big({\bf z}, \bl, \bb \Big) \rightarrow \Big(\tilde{\bf z}, \bl, \tilde{\bb} \Big)$$
using the transformation
$$\tilde{z_i} = z_i - {\bf x}_i^T\bb_a \ , \ \text{for all} \ i=1,2,\ldots, n \quad \text{and} \quad \tilde{\bb} = \bb - \bb_a$$
The absolute value of the Jacobian of this transformation is one, and the joint posterior density of $\left(\tilde{\bb}, \bl, 
\tilde{\bf z}\right)$ is given by 
\begin{align}
    &\pi\left(\tilde{\bb}, \bl, \tilde{\bf z} \Big| \y \right)\\
    &= \frac{1}{m(\y)}\Bigg[\prod_{i=1}^{n} \Big\{I_{\R_+}\left(\tilde{z}_i + {\bf x}_i^T\bb_a\right)I_{\{1\}}\left(y_i\right) + I_{\R_-}\left(\tilde{z}_i + {\bf x}_i^T\bb_a\right)I_{\{0\}}\left(y_i\right)\Big\}\\
    &\qquad\qquad \times \frac{\sqrt{\lambda_i}}{\sqrt{2\pi}}\exp\left\{-\frac{\lambda_i}{2}\left(\tilde{z}_i - {\bf x}_i^T\tilde{\bb}\right)^2\right\} \times \frac{\left(\frac{\nu}{2}\right)^{\frac{\nu}{2}}}{\Gamma\left(\frac{\nu}{2}\right)} \ \lambda_i^{\frac{\nu}{2} - 1}\exp\left\{-\frac{\nu}{2}\lambda_i\right\}\Bigg]\\
    &\qquad\qquad \times \left(2\pi\right)^{-\frac{p}{2}} \sqrt{\text{det}\left(\Sigma_a\right)} \ \exp\left[-\frac{1}{2} \ \tilde{\bb}^T\Sigma_a\tilde{\bb}\right]
    \label{eq:joint.density.BetaTilde_Lambda_ZTilde.TraceClass}
\end{align}

\noindent
Straightforward calculations using \eqref{eq:joint.density.BetaTilde_Lambda_ZTilde.TraceClass} show that 
\begin{equation}
    \tilde{\bb} \Big| \bl, \tilde{\bf z}, \y \ \sim \ \mathcal{N}_p\left(\left(X^T \Lambda X + \Sigma_a\right)^{-1}X^T\Lambda\tilde{\bf z} \ , \ \left(X^T \Lambda X + \Sigma_a\right)^{-1}\right)
    \label{eq:conditional.distribution.BetaTilde.Given.else.TraceClass}
\end{equation}

\noindent
and 
\begin{align}
    \pi\left(\tilde{\bb} \Big| \bl, \tilde{\bf z}, \y \right) &= \left(2\pi\right)^{-\frac{p}{2}} \sqrt{\text{det}\left(X^T \Lambda X + \Sigma_a\right)} \ \exp\Bigg[-\frac{1}{2}\Bigg\{\tilde{\bb}^T\left(X^T \Lambda X + \Sigma_a\right)\tilde{\bb} - 2\tilde{\bb}^TX^T\Lambda \tilde{\bf z}\\
    &\qquad\qquad\qquad\qquad\qquad\qquad\qquad + \tilde{\bf z}^T\Lambda X \left(X^T \Lambda X + \Sigma_a\right)^{-1}X^T\Lambda \tilde{\bf z}\Bigg\}\Bigg]
    \label{eq:conditional.density.FullForm.BetaTilde.Given.else.TraceClass}
\end{align}

It is easy to notice from  \eqref{eq:joint.density.BetaTilde_Lambda_ZTilde.TraceClass} that $\left(\lambda_i, \tilde{z}_i\right)$'s are conditionally independent given $\left(\tilde{\bb}, \y\right)$, and moreover,
\begin{align}
    \pi\left(\lambda_i, \tilde{z}_i \Big| \tilde{\bb}, y_i\right) &\propto \Big\{I_{\R_+}\left(\tilde{z}_i + {\bf x}_i^T\bb_a\right)I_{\{1\}}\left(y_i\right) + I_{\R_-}\left(\tilde{z}_i + {\bf x}_i^T\bb_a\right)I_{\{0\}}\left(y_i\right)\Big\}\\
    &\qquad\times\lambda_i^{\frac{\nu+1}{2} - 1} \ \exp\left[-\frac{\lambda_i}{2}\left\{\nu + \left(\tilde{z}_i - {\bf x}_i^T\tilde{\bb}\right)^2\right\}\right]. 
    \label{eq:conditional.density.lambda_i.z_iTilde.Given.else.TraceClass}
\end{align}

\noindent
Hence, $\lambda_1, \lambda_2, \ldots, \lambda_n$ are 
conditionally independent given $\left(\tilde{\bf z}, \tilde{\bb}, \y\right)$, and 
\begin{equation}
    \lambda_i \Big| \tilde{z}_i, \tilde{\bb}, y_i \ \sim \ \text{Gamma}\left(\frac{\nu+1}{2} \ , \ \frac{\nu + \left(\tilde{z}_i - {\bf x}_i^T\tilde{\bb}\right)^2}{2}\right)
    \label{eq:eq:conditional.distribution.lambda_i.Given.else.TraceClass}
\end{equation}

\noindent
for each $i\in\{1,2,\ldots,n\}$, which implies that 
\begin{align}
    \pi\left(\bl  \Big| \tilde{\bf z}, \tilde{\bb}, \y\right) &= K_1\left[\prod_{i=1}^{n}\left(\frac{\nu + \left(\tilde{z}_i - {\bf x}_i^T\tilde{\bb}\right)^2}{2}\right)^{\frac{\nu + 1}{2}}\right]\times \left[\prod_{i=1}^{n}\lambda_i^{\frac{\nu - 1}{2}}\right]\\
    &\qquad\qquad \times \exp\left[-\frac{1}{2}\sum_{i=1}^{n}\lambda_i\left(\nu + \left(\tilde{z}_i - {\bf x}_i^T\tilde{\bb}\right)^2\right)\right]\\
    &= K_1^{'}\left[\prod_{i=1}^{n}\left(\nu + \left(\tilde{z}_i - {\bf x}_i^T\tilde{\bb}\right)^2\right)^{\frac{\nu + 1}{2}}\right]\times \left[\prod_{i=1}^{n}\lambda_i^{\frac{\nu - 1}{2}}\right]\times \  \exp\left[-\frac{\nu}{2}\sum_{i=1}^{n}\lambda_i\right]\\
    &\qquad\qquad \times \exp\left[-\frac{1}{2}\left\{ \tilde{\bf z}^T\Lambda \tilde{\bf z} - 2\tilde{\bb}^TX^T\Lambda \tilde{\bf z} + \tilde{\bb}^TX^T \Lambda X\tilde{\bb}\right\}\right]
    \label{eq:conditional.density.FullForm.lambda.Given.else.TraceClass}
\end{align}

\noindent
where $K_1$ and $K_1'$ are appropriate constants. 

To find $\pi\left(\tilde{\bf z} \Big| \tilde{\bb}, \y\right)$, we use \eqref{eq:conditional.density.lambda_i.z_iTilde.Given.else.TraceClass} to get for each $i\in\{1,2,\ldots,n\}$,
\begin{align}
    &\pi\left(\tilde{z}_i \Big| \tilde{\bb}, y_i\right)\\
    &= \int_{\R_+}\pi\left(\lambda_i, \tilde{z}_i \Big| \tilde{\bb}, y_i\right) d\lambda_i\\
    &= C_{\tilde{\bb}, y_i} \times \left[I_{\R_+}\left(\tilde{z}_i + {\bf x}_i^T\bb_a\right)I_{\{1\}}\left(y_i\right) + I_{\R_-}\left(\tilde{z}_i + {\bf x}_i^T\bb_a\right)I_{\{0\}}\left(y_i\right)\right]\\
    &\qquad \times \int_{\R_+}\lambda_i^{\frac{\nu+1}{2} - 1} \ \exp\left[-\frac{\lambda_i}{2}\left\{\nu + \left(\tilde{z}_i - {\bf x}_i^T\tilde{\bb}\right)^2\right\}\right] d\lambda_i\\
    &\left[\text{where, }C_{\tilde{\bb}, y_i} \ \text{is a constant which depends on }\tilde{\bb}, y_i\right]\\
    &= C_{\tilde{\bb}, y_i} \times \left[I_{\R_+}\left(\tilde{z}_i + {\bf x}_i^T\bb_a\right)I_{\{1\}}\left(y_i\right) + I_{\R_-}\left(\tilde{z}_i + {\bf x}_i^T\bb_a\right)I_{\{0\}}\left(y_i\right)\right] \times \frac{\Gamma\left(\frac{\nu+1}{2}\right)}{\left(\frac{\nu + \left(\tilde{z}_i - {\bf x}_i^T\tilde{\bb}\right)^2}{2}\right)^{\frac{\nu+1}{2}}}\\
    &= C_{\tilde{\bb}, y_i}^{'} \times \left[I_{\R_+}\left(\tilde{z}_i + {\bf x}_i^T\bb_a\right)I_{\{1\}}\left(y_i\right) + I_{\R_-}\left(\tilde{z}_i + {\bf x}_i^T\bb_a\right)I_{\{0\}}\left(y_i\right)\right] \times \left(1+\frac{\left(\tilde{z}_i - {\bf x}_i^T\tilde{\bb}\right)^2}{\nu}\right)^{-\frac{\nu+1}{2}}
    \label{eq:conditional.density.z_iTilde.Given.betaTilde.y_i.TraceClass}
\end{align}

\noindent
where $C_{\tilde{\bb}, y_i}^{'}$ is the product of all constant terms that are free of $\tilde{z}_i$. We conclude from \eqref{eq:joint.density.BetaTilde_Lambda_ZTilde.TraceClass} and \eqref{eq:conditional.density.z_iTilde.Given.betaTilde.y_i.TraceClass} that conditional on $\left(\tilde{\bb}, \y\right)$, $\tilde{z}_1, \tilde{z}_2, \ldots, \tilde{z}_n$ are independent with $\tilde{z}_i \Big| \tilde{\bb}, \y$ following a truncated $t$ distribution with location ${\bf x}_i^T\tilde{\bb}$, scale $1$ and degrees of freedom $\nu$ that is truncated left at $-{\bf x}_i^T\bb_a$ if $y_i=1$ and truncated right at $-{\bf x}_i^T\bb_a$ if $y_i=0$. 

Now, if we denote $t_{\nu}\left(\mu, 1\right)$ to be the univariate Student's $t$ distribution with location $\mu$, scale $1$ and degrees of freedom $\nu$, and $F_{\nu}$ to be the cdf of the $t_{\nu}(0,1)$ distribution, then for $y_i = 0$,
\begin{align}
    \pi\left(\tilde{z}_i \Big| \tilde{\bb}, y_i\right) &= \frac{K_2\left(1+\frac{\left(\tilde{z}_i - {\bf x}_i^T\tilde{\bb}\right)^2}{\nu}\right)^{-\frac{\nu + 1}{2}}}{P\left(t_{\nu}\left({\bf x}_i^T\tilde{\bb}, \ 1\right)\leq -{\bf x}_i^T\bb_a\right)}\\
    &= \frac{K_2\left(1+\frac{\left(\tilde{z}_i - {\bf x}_i^T\tilde{\bb}\right)^2}{\nu}\right)^{-\frac{\nu + 1}{2}}}{P\left(t_{\nu}\left(0, \ 1\right)\leq -{\bf x}_i^T\left(\tilde{\bb}+\bb_a\right)\right)}\\
    &= \frac{K_2\left(1+\frac{\left(\tilde{z}_i - {\bf x}_i^T\tilde{\bb}\right)^2}{\nu}\right)^{-\frac{\nu + 1}{2}}}{1 - F_{\nu}\left({\bf x}_i^T\left(\tilde{\bb}+\bb_a\right)\right)}
    \label{eq:conditional.density.z_iTilde.Given.betaTilde.y_i=0.TraceClass}
\end{align}
and, for $y_i=1$,
\begin{align}
    \pi\left(\tilde{z}_i \Big| \tilde{\bb}, y_i\right) &= \frac{K_2\left(1+\frac{\left(\tilde{z}_i - {\bf x}_i^T\tilde{\bb}\right)^2}{\nu}\right)^{-\frac{\nu + 1}{2}}}{P\left(t_{\nu}\left({\bf x}_i^T\tilde{\bb}, \ 1\right)\geq -{\bf x}_i^T\bb_a\right)}\\
    &= \frac{K_2\left(1+\frac{\left(\tilde{z}_i - {\bf x}_i^T\tilde{\bb}\right)^2}{\nu}\right)^{-\frac{\nu + 1}{2}}}{P\left(t_{\nu}\left(0, \ 1\right)\geq -{\bf x}_i^T\left(\tilde{\bb}+\bb_a\right)\right)}\\
    &= \frac{K_2\left(1+\frac{\left(\tilde{z}_i - {\bf x}_i^T\tilde{\bb}\right)^2}{\nu}\right)^{-\frac{\nu + 1}{2}}}{F_{\nu}\left({\bf x}_i^T\left(\tilde{\bb}+\bb_a\right)\right)}. 
    \label{eq:conditional.density.z_iTilde.Given.betaTilde.y_i=1.TraceClass}
\end{align}
Combining \eqref{eq:conditional.density.z_iTilde.Given.betaTilde.y_i=0.TraceClass} and \eqref{eq:conditional.density.z_iTilde.Given.betaTilde.y_i=1.TraceClass}, we have for any $y_i$,
\begin{align}
   &\pi\left(\tilde{z}_i \Big| \tilde{\bb}, y_i\right)\\
   &= K_2\left(1+\frac{\left(\tilde{z}_i - {\bf x}_i^T\tilde{\bb}\right)^2}{\nu}\right)^{-\frac{\nu + 1}{2}} \left\{\frac{1}{F_{\nu}\left({\bf x}_i^T\left(\tilde{\bb}+\bb_a\right)\right)}\right\}^{y_i}\left\{\frac{1}{1 - F_{\nu}\left({\bf x}_i^T\left(\tilde{\bb}+\bb_a\right)\right)}\right\}^{1-y_i}
\end{align}
which implies,
\begin{align}
    &\pi\left(\tilde{\bf z} \Big| \tilde{\bb}, \y\right)\\
    &= K_2^n \ \nu^{\frac{n\left(\nu+1\right)}{2}} \ \prod_{i=1}^{n}\Bigg[\left(\nu + \left(\tilde{z}_i - {\bf x}_i^T\tilde{\bb}\right)^2\right)^{-\frac{\nu + 1}{2}}\left\{\frac{1}{F_{\nu}\left({\bf x}_i^T\left(\tilde{\bb}+\bb_a\right)\right)}\right\}^{y_i}\\
    &\qquad\qquad\qquad\qquad\qquad\qquad\times \left\{\frac{1}{1 - F_{\nu}\left({\bf x}_i^T\left(\tilde{\bb}+\bb_a\right)\right)}\right\}^{1-y_i}\Bigg]. 
    \label{eq:conditional.density.FullForm.zTilde.Given.betaTilde.y.TraceClass}
\end{align}

\noindent
Let $S:= \left\{i: y_i = 0\right\}$. Then, $S^c = \left\{i : \ y_i = 1\right\}$. Using \eqref{eq:markov.transition.density}, \eqref{eq:conditional.density.FullForm.BetaTilde.Given.else.TraceClass}, \eqref{eq:conditional.density.FullForm.lambda.Given.else.TraceClass} and \eqref{eq:conditional.density.FullForm.zTilde.Given.betaTilde.y.TraceClass}, we get the following form for the integral $I$ in \eqref{eq:trace.class.condition} under 
the new parametrization. 
\begin{align}
    I &:= \int_{\R^p}k\left(\tilde{\bb}, \tilde{\bb}\right)d\tilde{\bb}\\
    &= \int_{\R^p}\int_{\R^n_{+}}\int_{\R^n}\pi\left(\tilde{\bb} \Big| \bl, \tilde{\bf z}, \y\right)\pi\left(\bl \Big| \tilde{\bf z}, \tilde{\bb}, \y\right)\pi\left(\tilde{\bf z} \Big| \tilde{\bb}, \y\right)d\tilde{\bf z} \ d\bl \ d\tilde{\bb}\\
    &= C_0 \int_{\R^p}\int_{\R^n_{+}}\int_{\R^n}\sqrt{\text{det}\left(X^T \Lambda X + \Sigma_a\right)} \ \times \left[\prod_{i=1}^{n}\lambda_i^{\frac{\nu - 1}{2}}\right]\times \  \exp\left[-\frac{\nu}{2}\sum_{i=1}^{n}\lambda_i\right]\\
    &\qquad\qquad\qquad\qquad \times \prod_{i\in S}\left[\frac{1}{1 - F_{\nu}\left({\bf x}_i^T\left(\tilde{\bb}+\bb_a\right)\right)}\right] \times \prod_{i\in S^c}\left[\frac{1}{F_{\nu}\left({\bf x}_i^T\left(\tilde{\bb}+\bb_a\right)\right)}\right]\\
    &\qquad\qquad\qquad\qquad \times \ \exp\Bigg[-\frac{1}{2}\Bigg\{\tilde{\bb}^T\left(X^T \Lambda X + \Sigma_a\right)\tilde{\bb} \ - \ 2\tilde{\bb}^TX^T\Lambda \tilde{\bf z}\\
    &\qquad\qquad\qquad\qquad\qquad\qquad\qquad + \ \tilde{\bf z}^T\Lambda X \left(X^T \Lambda X + \Sigma_a\right)^{-1}X^T\Lambda \tilde{\bf z}\Bigg\}\Bigg]\\
    &\qquad\qquad\qquad\qquad \times \ \exp\left[-\frac{1}{2}\left\{ \tilde{\bf z}^T\Lambda \tilde{\bf z} - 2\tilde{\bb}^TX^T\Lambda \tilde{\bf z} + \tilde{\bb}^TX^T \Lambda X\tilde{\bb}\right\}\right] \ d\tilde{\bf z} \ d\bl \ d\tilde{\bb}\\
    \label{eq:integral.value.1.TraceClass}
\end{align}
Here, $C_0$ denotes the product of all constant terms (independent of $\tilde{\bb}$, $\bl$, and $\tilde{\bf z}$) appearing in the conditional densities $\pi\left(\tilde{\bb} \Big| \bl, \tilde{\bf z}, \y\right)$, $\pi\left(\bl \Big| \tilde{\bf z}, \tilde{\bb}, \y \right)$, and $\pi\left(\tilde{\bf z} \Big| \tilde{\bb}, \y\right)$.

\medskip

\noindent {\bf Step II: Another reparametrization to adjust for 
the prior precision matrix $\Sigma_a$.} Now, let us define ${\boldsymbol \theta} = \Sigma_a^{1/2}\tilde{\bb}$, $W = X\Sigma_a^{-1/2}$, and $\tilde{\bf c} = \Sigma_a^{1/2}\bb_a$. Absolute value of the Jacobian of the transformation \ $\tilde{\bb}\rightarrow {\boldsymbol \theta}$ \ is $\left\{\text{det}\left(\Sigma_a\right)\right\}^{-1/2} > 0$. Therefore, the right hand side of \eqref{eq:integral.value.1.TraceClass}, after this 
further reparametrization becomes
\begin{align}
    I &= C_0 \int_{\R^p}\int_{\R^n_{+}}\int_{\R^n}\sqrt{\text{det}\left(W^T\Lambda W + I_p\right)} \ \times \left[\prod_{i=1}^{n}\lambda_i^{\frac{\nu - 1}{2}}\right]\times \  \exp\left[-\frac{\nu}{2}\sum_{i=1}^{n}\lambda_i\right]\\
    &\qquad\qquad\qquad\qquad \times \prod_{i\in S}\left[\frac{1}{1 - F_{\nu}\left({\bf w}_i^T\left({\boldsymbol \theta}+\tilde{\bf c}\right)\right)}\right] \times \prod_{i\in S^c}\left[\frac{1}{F_{\nu}\left({\bf w}_i^T\left({\boldsymbol \theta}+\tilde{\bf c}\right)\right)}\right]\\
    &\qquad\qquad\qquad\qquad \times \ \exp\Bigg[-\frac{1}{2}\Bigg\{{\boldsymbol \theta}^T\left(W^T\Lambda W + I_p\right){\boldsymbol \theta} \ - \ 2{\boldsymbol \theta}^TW^T\Lambda \tilde{\bf z}\\
    &\qquad\qquad\qquad\qquad\qquad\qquad\qquad + \ \tilde{\bf z}^T\Lambda  W \left(W^T\Lambda W + I_p\right)^{-1}W^T\Lambda \tilde{\bf z}\Bigg\}\Bigg]\\
    &\qquad\qquad\qquad\qquad \times \ \exp\left[-\frac{1}{2}\left\{ \tilde{\bf z}^T\Lambda \tilde{\bf z} - 2{\boldsymbol \theta}^TW^T\Lambda \tilde{\bf z} + {\boldsymbol \theta}^TW^T\Lambda W{\boldsymbol \theta}\right\}\right] \ d\tilde{\bf z} \ d\bl \ d{\boldsymbol \theta}\\
    &= C_0 \int_{\R^p}\int_{\R^n_{+}}\sqrt{\text{det}\left(W^T\Lambda W + I_p\right)} \ \times \left[\prod_{i=1}^{n}\lambda_i^{\frac{\nu - 1}{2}}\right]\times \  \exp\left[-\frac{\nu}{2}\sum_{i=1}^{n}\lambda_i\right]\\
    &\qquad\quad \times \prod_{i\in S}\left[\frac{1}{1 - F_{\nu}\left({\bf w}_i^T\left({\boldsymbol \theta}+\tilde{\bf c}\right)\right)}\right] \times \prod_{i\in S^c}\left[\frac{1}{F_{\nu}\left({\bf w}_i^T\left({\boldsymbol \theta}+\tilde{\bf c}\right)\right)}\right]\\
    &\qquad\quad \times \ \exp\Bigg[-\frac{1}{2}\Bigg\{{\boldsymbol \theta}^T\left(2W^T\Lambda W + I_p\right){\boldsymbol \theta}\Bigg\}\Bigg]\\
    &\qquad\quad \times \left(\int_{\R^n}\exp\left[2{\boldsymbol \theta}^TW^T\Lambda \tilde{\bf z} - \frac{1}{2} \ \tilde{\bf z}^T\left(\Lambda + \Lambda  W \left(W^T\Lambda W + I_p\right)^{-1}W^T\Lambda\right) \tilde{\bf z}\right] d\tilde{\bf z}\right)\\
    &\qquad\qquad\qquad\qquad\qquad\qquad\qquad\qquad\qquad\qquad\qquad\qquad\qquad\qquad\qquad\qquad d\bl \ d{\boldsymbol \theta}\\
    \label{eq:integral.value.2.TraceClass}
\end{align}

\noindent
{\bf Step III: An upper bound for the innermost $\tilde{\bf z}$ integral in
\eqref{eq:integral.value.2.TraceClass}.} We now derive an upper bound
for the innermost integral in \eqref{eq:integral.value.2.TraceClass}. Note 
that 
\begin{align}
    & \int_{\R^n}\exp\left[2{\boldsymbol \theta}^TW^T\Lambda \tilde{\bf z} \ - \ \frac{1}{2} \ \tilde{\bf z}^T\left(\Lambda + \Lambda  W \left(W^T\Lambda W + I_p\right)^{-1}W^T\Lambda\right) \tilde{\bf z}\right] \ d\tilde{\bf z}\\
    &= \ \exp\left[\frac{1}{2} \ 4 {\boldsymbol \theta}^TW^T\Lambda\left(\Lambda + \Lambda  W \left(W^T\Lambda W + I_p\right)^{-1}W^T\Lambda\right)^{-1}\Lambda  W{\boldsymbol \theta}\right] \times \left(C_1\right)^{-1}\\
    &\qquad\qquad \times \int_{\R^n}C_1 \exp\left[-\frac{1}{2}\left(\tilde{\bf z} - {\bf a}_1^{\star}\right)^T\left(\Lambda + \Lambda  W \left(W^T\Lambda W + I_p\right)^{-1}W^T\Lambda\right)\left(\tilde{\bf z} - {\bf a}_1^{\star}\right)\right] \ d\tilde{\bf z}\\
    &\qquad\qquad \left[\text{where, }{\bf a}_1^{\star} = 2\left(\Lambda + \Lambda  W \left(W^T\Lambda W + I_p\right)^{-1}W^T\Lambda\right)^{-1}\Lambda  W{\boldsymbol \theta}\right]\\
    &= \ \left(C_1\right)^{-1} \exp\left[\frac{1}{2} \ 4 {\boldsymbol \theta}^TW^T\Lambda\left(\Lambda + \Lambda  W \left(W^T\Lambda W + I_p\right)^{-1}W^T\Lambda\right)^{-1}\Lambda  W{\boldsymbol \theta}\right]
    \label{eq:inner.integral.value.TraceClass}
\end{align}
where, $C_1 = \left(2\pi\right)^{-n/2}\left\{\text{det}\left(\Lambda + \Lambda  W \left(W^T\Lambda W + I_p\right)^{-1}W^T\Lambda\right)\right\}^{1/2}$. The last equality follows from the fact that the integrand is a normal density. However,
\begin{align}
    C_1 &= \left(2\pi\right)^{-n/2}\left\{\text{det}\left(\Lambda + \Lambda  W \left(W^T\Lambda W + I_p\right)^{-1}W^T\Lambda\right)\right\}^{1/2}\\
    &\geq \left(2\pi\right)^{-n/2} \left\{\text{det}\left(\Lambda\right)\right\}^{1/2}\\
    &= \left(2\pi\right)^{-n/2} \sqrt{\prod_{i=1}^{n}\lambda_i}
    \label{eq:inner.integral.value.Normalizing.Constant.lower.bound.TraceClass}
\end{align}
From \eqref{eq:inner.integral.value.TraceClass} and \eqref{eq:inner.integral.value.Normalizing.Constant.lower.bound.TraceClass}, we have an upper bound for the inner integral in \eqref{eq:integral.value.2.TraceClass} as follows
\begin{align}
    & \int_{\R^n}\exp\left[2{\boldsymbol \theta}^TW^T\Lambda \tilde{\bf z} \ - \ \frac{1}{2} \ \tilde{\bf z}^T\left(\Lambda + \Lambda  W \left(W^T\Lambda W + I_p\right)^{-1}W^T\Lambda\right) \tilde{\bf z}\right] \ d\tilde{\bf z}\\
    &\leq \left(2\pi\right)^{n/2} \prod_{i=1}^{n}\lambda_i^{-1/2} \ \exp\left[\frac{1}{2} \ 4 {\boldsymbol \theta}^TW^T\Lambda\left(\Lambda + \Lambda  W \left(W^T\Lambda W + I_p\right)^{-1}W^T\Lambda\right)^{-1}\Lambda  W{\boldsymbol \theta}\right]\\
    & \label{eq:inner.integral.upper.bound.TraceClass}
\end{align}

\noindent
Combining \eqref{eq:integral.value.2.TraceClass} and \eqref{eq:inner.integral.upper.bound.TraceClass}, we 
get
\begin{align}
    I &\leq C_0\left(2\pi\right)^{n/2}\int_{\R^p}\int_{\R^n_{+}}\sqrt{\text{det}\left(W^T\Lambda W + I_p\right)} \ \times \left[\prod_{i=1}^{n}\lambda_i^{\frac{\nu}{2} - 1}\right]\times \  \exp\left[-\frac{\nu}{2}\sum_{i=1}^{n}\lambda_i\right]\\
    & \qquad\qquad\qquad \times \prod_{i\in S}\left[\frac{1}{1 - F_{\nu}\left({\bf w}_i^T\left({\boldsymbol \theta}+\tilde{\bf c}\right)\right)}\right] \times \prod_{i\in S^c}\left[\frac{1}{F_{\nu}\left({\bf w}_i^T\left({\boldsymbol \theta}+\tilde{\bf c}\right)\right)}\right]\\
    & \qquad\qquad\qquad \times \exp\Bigg[-\frac{1}{2}\Bigg\{{\boldsymbol \theta}^T\left(2W^T\Lambda W + I_p\right){\boldsymbol \theta}\\
    & \qquad\qquad\qquad\qquad\qquad - 4 {\boldsymbol \theta}^TW^T\Lambda\left(\Lambda + \Lambda  W \left(W^T\Lambda W + I_p\right)^{-1}W^T\Lambda\right)^{-1}\Lambda  W{\boldsymbol \theta}\Bigg\}\Bigg]\\
    &\qquad\qquad\qquad\qquad\qquad\qquad\qquad\qquad\qquad\qquad\qquad\qquad\qquad\qquad\qquad\qquad\quad d\bl \ d{\boldsymbol \theta}\\
    &= C_2\int_{\R^p}\int_{\R^n_{+}}\sqrt{\text{det}\left(W^T\Lambda W + I_p\right)} \ \times \left[\prod_{i=1}^{n}\lambda_i^{\frac{\nu}{2} - 1}\right]\times \  \exp\left[-\frac{\nu}{2}\sum_{i=1}^{n}\lambda_i\right]\\
    & \qquad\qquad\qquad \times \prod_{i\in S}\left[\frac{1}{1 - F_{\nu}\left({\bf w}_i^T\left({\boldsymbol \theta}+\tilde{\bf c}\right)\right)}\right] \times \prod_{i\in S^c}\left[\frac{1}{F_{\nu}\left({\bf w}_i^T\left({\boldsymbol \theta}+\tilde{\bf c}\right)\right)}\right]\\
    & \qquad\qquad\qquad \times \exp\left[-\frac{1}{2}\Big\{G\left({\boldsymbol \theta}, \bl\right)\Big\}\right] \ d\bl \ d{\boldsymbol \theta}
    \label{eq:integral.value.3.TraceClass}
\end{align}
where $C_2 = C_0\left(2\pi\right)^{n/2}$ and 
\begin{align}
    G\left({\boldsymbol \theta}, \bl\right) &= {\boldsymbol \theta}^T\left(2W^T\Lambda W + I_p\right){\boldsymbol \theta} - 4 {\boldsymbol \theta}^TW^T\Lambda\left(\Lambda + \Lambda  W \left(W^T\Lambda W + I_p\right)^{-1}W^T\Lambda\right)^{-1}\Lambda  W{\boldsymbol \theta}\\
    &= {\boldsymbol \theta}^T\left[\left(2W^T\Lambda W + I_p\right) - 4 W^T\Lambda\left(\Lambda + \Lambda  W \left(W^T\Lambda W + I_p\right)^{-1}W^T\Lambda\right)^{-1}\Lambda  W\right]{\boldsymbol \theta}
    \label{eq:G.theta.lambda.function.TraceClass}
\end{align}

\noindent
{\bf Step IV: An upper bound for the products involving
the cdf $F_{\nu}$.} We now target the product terms in 
the integrand involving the $t$-cdf $F_\nu$. Note that 
for $i\in S$, if ${\bf w}_i^T\left({\boldsymbol \theta}+\tilde{\bf c}\right)\leq 0$, then
\begin{align}
    F_{\nu}\left({\bf w}_i^T\left({\boldsymbol \theta}+\tilde{\bf c}\right)\right) \leq F_{\nu}(0) = \frac{1}{2} 
    \implies \frac{1}{1 - F_{\nu}\left({\bf w}_i^T\left({\boldsymbol \theta}+\tilde{\bf c}\right)\right)} \leq 2. 
    \label{eq:corollary1.derivative.1}
\end{align}
and, if ${\bf w}_i^T\left({\boldsymbol \theta}+\tilde{\bf c}\right)> 0$, then by Corollary \ref{appendix:corollary.1} in Appendix \ref{appendix:mills.ratio.type.result} we have
\begin{align}
    \frac{1}{1 - F_{\nu}\left({\bf w}_i^T\left({\boldsymbol \theta}+\tilde{\bf c}\right)\right)} &\leq \frac{\left(\left({\bf w}_i^T{\boldsymbol \theta}+{\bf w}_i^T\tilde{\bf c}\right)^2+\nu\right)^{\frac{\nu}{2}}}{\kappa}. 
    \label{eq:corollary1.derivative.2}
\end{align}
From \eqref{eq:corollary1.derivative.1} and \eqref{eq:corollary1.derivative.2}, we have for any $i\in S$,
\begin{align}
    \frac{1}{1 - F_{\nu}\left({\bf w}_i^T\left({\boldsymbol \theta}+\tilde{\bf c}\right)\right)} &\leq \text{max}\left\{2 \ , \ \frac{\left(\left({\bf w}_i^T{\boldsymbol \theta}+{\bf w}_i^T\tilde{\bf c}\right)^2+\nu\right)^{\frac{\nu}{2}}}{\kappa}\right\}\\
    &\leq \left(2 + \frac{\left(\left({\bf w}_i^T{\boldsymbol \theta}+{\bf w}_i^T\tilde{\bf c}\right)^2+\nu\right)^{\frac{\nu}{2}}}{\kappa}\right). 
    \label{eq:corollary1.derivative.i.in.S}
\end{align}
Similarly for $i\in S^c$, if ${\bf w}_i^T\left({\boldsymbol \theta}+\tilde{\bf c}\right)\geq 0$, then
\begin{align}
    F_{\nu}\left({\bf w}_i^T\left({\boldsymbol \theta}+\tilde{\bf c}\right)\right) &\geq F_{\nu}(0) = \frac{1}{2}\\
    \implies \frac{1}{F_{\nu}\left({\bf w}_i^T\left({\boldsymbol \theta}+\tilde{\bf c}\right)\right)} &\leq 2. 
    \label{eq:corollary1.derivative.3}
\end{align}
and, if ${\bf w}_i^T\left({\boldsymbol \theta}+\tilde{\bf c}\right)< 0$, i.e., $-{\bf w}_i^T\left({\boldsymbol \theta}+\tilde{\bf c}\right)> 0$, then by Corollary \ref{appendix:corollary.1} in 
Appendix \ref{appendix:mills.ratio.type.result} we have
\begin{equation}
    \frac{1}{F_{\nu}\left({\bf w}_i^T\left({\boldsymbol \theta}+\tilde{\bf c}\right)\right)} \ = \ \frac{1}{1-F_{\nu}\left(-{\bf w}_i^T\left({\boldsymbol \theta}+\tilde{\bf c}\right)\right)} \ \leq \ \frac{\left(\left({\bf w}_i^T{\boldsymbol \theta}+{\bf w}_i^T\tilde{\bf c}\right)^2+\nu\right)^{\frac{\nu}{2}}}{\kappa}
    \label{eq:corollary1.derivative.4}
\end{equation}
From \eqref{eq:corollary1.derivative.3} and \eqref{eq:corollary1.derivative.4}, we have for any $i\in S^c$,
\begin{align}
    \frac{1}{F_{\nu}\left({\bf w}_i^T\left({\boldsymbol \theta}+\tilde{\bf c}\right)\right)} &\leq \text{max}\left\{2 \ , \ \frac{\left(\left({\bf w}_i^T{\boldsymbol \theta}+{\bf w}_i^T\tilde{\bf c}\right)^2+\nu\right)^{\frac{\nu}{2}}}{\kappa}\right\}\\
    &\leq \left(2 + \frac{\left(\left({\bf w}_i^T{\boldsymbol \theta}+{\bf w}_i^T\tilde{\bf c}\right)^2+\nu\right)^{\frac{\nu}{2}}}{\kappa}\right)
    \label{eq:corollary1.derivative.i.in.S.complement}
\end{align}
Finally from \eqref{eq:corollary1.derivative.i.in.S} and \eqref{eq:corollary1.derivative.i.in.S.complement}, we have
\begin{align}
    \prod_{i\in S}\left[\frac{1}{1 - F_{\nu}\left({\bf w}_i^T\left({\boldsymbol \theta}+\tilde{\bf c}\right)\right)}\right] \times \prod_{i\in S^c}\left[\frac{1}{F_{\nu}\left({\bf w}_i^T\left({\boldsymbol \theta}+\tilde{\bf c}\right)\right)}\right] &\leq \prod_{i=1}^{n}\left(2 + \frac{\left(\left({\bf w}_i^T{\boldsymbol \theta}+{\bf w}_i^T\tilde{\bf c}\right)^2+\nu\right)^{\frac{\nu}{2}}}{\kappa}\right)
    \label{eq:corollary1.derivative.upper.bound.1.TraceClass}
\end{align}

\noindent
Note that 
\begin{align}
    \left({\bf w}_i^T{\boldsymbol \theta}+{\bf w}_i^T\tilde{\bf c}\right)^2 &= \left({\boldsymbol \theta}+\tilde{\bf c}\right)^T{\bf w}_i{\bf w}_i^T\left({\boldsymbol \theta}+\tilde{\bf c}\right)\\
    &\leq \left({\boldsymbol \theta}+\tilde{\bf c}\right)^T\left(\sum_{i=1}^{n}{\bf w}_i{\bf w}_i^T\right)\left({\boldsymbol \theta}+\tilde{\bf c}\right)\\
    &= \left({\boldsymbol \theta}+\tilde{\bf c}\right)^TW^TW\left({\boldsymbol \theta}+\tilde{\bf c}\right)
    \label{eq:observation.1.TraceClass}
\end{align}
for every $1 \leq i \leq n$. It follows from  \eqref{eq:corollary1.derivative.upper.bound.1.TraceClass}, \eqref{eq:observation.1.TraceClass}, and the 
$c_r$-inequality that 
\begin{align}
    & \prod_{i\in S}\left[\frac{1}{1 - F_{\nu}\left({\bf w}_i^T\left({\boldsymbol \theta}+\tilde{\bf c}\right)\right)}\right] \times \prod_{i\in S^c}\left[\frac{1}{F_{\nu}\left({\bf w}_i^T\left({\boldsymbol \theta}+\tilde{\bf c}\right)\right)}\right] \\ &\leq \ \prod_{i=1}^{n}\left(2 + \frac{\left(\left({\boldsymbol \theta}+\tilde{\bf c}\right)^TW^TW\left({\boldsymbol \theta}+\tilde{\bf c}\right)+\nu\right)^{\frac{\nu}{2}}}{\kappa}\right)\\
    &= \ \left(2 + \frac{\left(\left({\boldsymbol \theta}+\tilde{\bf c}\right)^TW^TW\left({\boldsymbol \theta}+\tilde{\bf c}\right)+\nu\right)^{\frac{\nu}{2}}}{\kappa}\right)^n\\
    &\leq \ 2^n \left[2^n \ + \ \frac{\left(\left({\boldsymbol \theta}+\tilde{\bf c}\right)^TW^TW\left({\boldsymbol \theta}+\tilde{\bf c}\right)+\nu\right)^{\frac{n\nu}{2}}}{\kappa^n}\right]\\
    &\leq \ 2^n \left[2^n \ + \ \frac{2^{n\nu/2}\left\{\left(\left({\boldsymbol \theta}+\tilde{\bf c}\right)^TW^TW\left({\boldsymbol \theta}+\tilde{\bf c}\right)\right)^{\frac{n\nu}{2}}+\nu^{n\nu/2}\right\}}{\kappa^n}\right]\\
    &= \ C_3 \ + \ C_4\left(\left({\boldsymbol \theta}+\tilde{\bf c}\right)^TW^TW\left({\boldsymbol \theta}+\tilde{\bf c}\right)\right)^{\frac{n\nu}{2}}\\
    &= \ C_3 \ + \ C_6\left(\left({\boldsymbol \theta}+\tilde{\bf c}\right)^T\left({\boldsymbol \theta}+\tilde{\bf c}\right)\right)^{\frac{n\nu}{2}}\label{eq:corollary1.derivative.upper.bound.3.TraceClass}, 
\end{align}

\noindent
where $C_3 = 2^{2n} + 
\frac{2^n\left(2\nu\right)^{n\nu/2}}{\kappa^n}$, 
$C_4 = \frac{2^n2^{n\nu/2}}{\kappa^n}$, $C_5$ denotes  
the largest eigenvalue of $W^T W$, and $C_6 = 
C_4C_5^{\frac{n\nu}{2}}$. 

Now, from \eqref{eq:integral.value.3.TraceClass},  \eqref{eq:corollary1.derivative.upper.bound.3.TraceClass} and Fubini's theorem, we get
\begin{align}
    I \ &\leq \ C_2\int_{\R^p}\int_{\R^n_{+}}\sqrt{\text{det}\left(W^T\Lambda W + I_p\right)} \ \times \left[\prod_{i=1}^{n}\lambda_i^{\frac{\nu}{2} - 1}\right]\times \  \exp\left[-\frac{\nu}{2}\sum_{i=1}^{n}\lambda_i\right]\\
    & \qquad\qquad\qquad\qquad \times \left[C_3 \ + \ C_6\left(\left({\boldsymbol \theta}+\tilde{\bf c}\right)^T\left({\boldsymbol \theta}+\tilde{\bf c}\right)\right)^{\frac{n\nu}{2}}\right]\\
    & \qquad\qquad\qquad\qquad \times \exp\left[-\frac{1}{2}\Big\{G\left({\boldsymbol \theta}, \bl\right)\Big\}\right] \ d\bl \ d{\boldsymbol \theta}\\
    &= \ C_2\int_{\R^n_{+}}\int_{\R^p}\sqrt{\text{det}\left(W^T\Lambda W + I_p\right)} \ \times \left[\prod_{i=1}^{n}\lambda_i^{\frac{\nu}{2} - 1}\right]\times \  \exp\left[-\frac{\nu}{2}\sum_{i=1}^{n}\lambda_i\right]\\
    & \qquad\qquad\qquad\qquad \times \left[C_3 \ + \ C_6\left(\left({\boldsymbol \theta}+\tilde{\bf c}\right)^T\left({\boldsymbol \theta}+\tilde{\bf c}\right)\right)^{\frac{n\nu}{2}}\right]\\
    & \qquad\qquad\qquad\qquad \times \exp\left[-\frac{1}{2}\Big\{G\left({\boldsymbol \theta}, \bl\right)\Big\}\right] \ d{\boldsymbol \theta} \ d\bl
    \label{eq:integral.value.4.TraceClass}
\end{align}

\noindent
Let $\lambda_{\text{max}}$ denote the largest diagonal 
entry of the diagonal matrix $\Lambda$. Using 
$\left(W^T\Lambda W + I_p\right) \preceq 
\max(1,\lambda_{\text{max}})  \left(W^TW + 
I_p\right)$, it follows that
\begin{align}
    I \ &\leq \ C_2\int_{\R^n_{+}}\int_{\R^p}\text{max}\left\{1,\lambda_{\text{max}}^{p/2}\right\} \sqrt{\text{det}\left(W^TW + I_p\right)} \ \times \left[\prod_{i=1}^{n}\lambda_i^{\frac{\nu}{2} - 1}\right]\times \  \exp\left[-\frac{\nu}{2}\sum_{i=1}^{n}\lambda_i\right]\\
    & \qquad\qquad\qquad\qquad \times \left[C_3 \ + \ C_6\left(\left({\boldsymbol \theta}+\tilde{\bf c}\right)^T\left({\boldsymbol \theta}+\tilde{\bf c}\right)\right)^{\frac{n\nu}{2}}\right]\\
    & \qquad\qquad\qquad\qquad \times \exp\left[-\frac{1}{2}\Big\{G\left({\boldsymbol \theta}, \bl\right)\Big\}\right] \ d{\boldsymbol \theta} \ d\bl\\
    &= \ C_7\int_{\R^n_{+}}\int_{\R^p}\text{max}\left\{1,\lambda_{\text{max}}^{p/2}\right\} \times \left[\prod_{i=1}^{n}\lambda_i^{\frac{\nu}{2} - 1}\right]\times \  \exp\left[-\frac{\nu}{2}\sum_{i=1}^{n}\lambda_i\right]\\
    & \qquad\qquad\qquad\qquad \times \left[C_3 \ + \ C_6\left(\left({\boldsymbol \theta}+\tilde{\bf c}\right)^T\left({\boldsymbol \theta}+\tilde{\bf c}\right)\right)^{\frac{n\nu}{2}}\right]\\
    & \qquad\qquad\qquad\qquad \times \exp\left[-\frac{1}{2}\Big\{G\left({\boldsymbol \theta}, \bl\right)\Big\}\right] \ d{\boldsymbol \theta} \ d\bl
    \label{eq:integral.value.5.TraceClass}
\end{align}

\noindent
where $C_7 = C_2\sqrt{\text{det}\left(W^TW + 
I_p\right)}$ is a constant term free of ${\boldsymbol 
\theta}$ and $\bl$. 

\iffalse
\begin{align}
    I_1 \ &= \ C_3C_7\int_{\R^n_{+}}\int_{\R^p}\text{max}\left\{1,\lambda_{\text{max}}^{p/2}\right\} \times \left[\prod_{i=1}^{n}\lambda_i^{\frac{\nu}{2} - 1}\right]\times \  \exp\left[-\frac{\nu}{2}\sum_{i=1}^{n}\lambda_i\right]\\
    & \qquad\qquad\qquad\qquad \times \exp\Big[-\frac{1}{2}\Big\{G\left({\boldsymbol \theta}, \lambda\right)\Big\}\Big] \ d{\boldsymbol \theta} \ d\bl
    \label{eq:integral.value.5_1.TraceClass}
\end{align}
and
\begin{align}
    I_2 \ &= \ C_6C_7\int_{\R^n_{+}}\int_{\R^p}\text{max}\left\{1,\lambda_{\text{max}}^{p/2}\right\} \times \left[\prod_{i=1}^{n}\lambda_i^{\frac{\nu}{2} - 1}\right]\times \  \exp\left[-\frac{\nu}{2}\sum_{i=1}^{n}\lambda_i\right]\\
    & \qquad\qquad\qquad\qquad \times \left(\left({\boldsymbol \theta}+\tilde{c}\right)^T\left({\boldsymbol \theta}+\tilde{c}\right)\right)^{\frac{n\nu}{2}} \times \exp\Big[-\frac{1}{2}\Big\{G\left({\boldsymbol \theta}, \lambda\right)\Big\}\Big] \ d{\boldsymbol \theta} \ d\bl. 
    \label{eq:integral.value.5_2.TraceClass}
\end{align}
Therefore to prove \eqref{eq:trace.class.condition}, it is enough to show that both $I_1$ and $I_2$ are finite.
\fi

\medskip

\noindent
{\bf Step V: Showing $G\left({\boldsymbol \theta}, \bl\right)$ is positive definite quadratic form in 
${\boldsymbol \theta}$.} In order to show the finiteness of the upper bound for $I$ in \eqref{eq:integral.value.5.TraceClass}, we will first prove that $G\left({\boldsymbol \theta}, \bl\right)$ is a positive definite quadratic form in ${\boldsymbol \theta}$ for all ${\boldsymbol \theta} \in \R^p$ and $\bl \in \R^n_{+}$. For that, it is enough to show by \eqref{eq:G.theta.lambda.function.TraceClass} that the matrix 
$$
\left[\left(2W^T\Lambda W + I_p\right) - 4 W^T\Lambda\left(\Lambda + \Lambda  W \left(W^T\Lambda W + I_p\right)^{-1}W^T\Lambda\right)^{-1}\Lambda  W\right]
$$

\noindent
is a positive definite matrix for all $\bl \in \R^n_{+}$. We show this by working out the spectral decomposition of this matrix separately in the low and high-dimensional settings. 

\medskip

\noindent
\textbf{\underline{Low-dimensional setting}:} When $n 
\geq p$
\begin{align}
    & \left(2W^T\Lambda W + I_p\right) - 4 W^T\Lambda\left(\Lambda + \Lambda  W \left(W^T\Lambda W + I_p\right)^{-1}W^T\Lambda\right)^{-1}\Lambda  W\\
    &= \left(2W^T\Lambda W + I_p\right) - 4 W^T\Lambda^{1/2}\left(I_n + \Lambda^{1/2} W \left(W^T\Lambda W + I_p\right)^{-1}W^T\Lambda^{1/2}\right)^{-1}\Lambda^{1/2} W\\
    &= \left(2A^TA + I_p\right) - 4 A^T\left(I_n + A \left(A^TA + I_p\right)^{-1}A^T\right)^{-1}A
    \label{eq:CaseI.1.TraceClass}
\end{align}
where $A = \Lambda^{1/2} W$. Now, by the Singular Value Decomposition, $A$ can be written as
\begin{equation}
    A_{n\times p} = U_{n\times p}D_{p\times p}V_{p\times p}^T
    \label{eq:SVD.A.CaseI.TraceClass}
\end{equation}
where, $U$ is a semi-orthogonal matrix i.e. $U^TU=I_p$, $V$ is an orthogonal matrix i.e. $V^TV = VV^T = I_p$, and $D$ is a diagonal matrix with singular values of $A$ being the diagonal entries. Since, $U_{n\times p}$ is a semi-orthogonal matrix with full column rank $p$, there exists a matrix $U_{{0}_{ \ n \times (n-p)}}$ such that the matrix $\begin{bmatrix}U_{n\times p} & U_{{0}_{ \ n \times (n-p)}} \end{bmatrix}$ becomes an orthogonal matrix, i.e.,
\begin{equation}
    \begin{bmatrix}U & U_0 \end{bmatrix}\begin{bmatrix}U^T\\ U_0^T \end{bmatrix} \ = \ \begin{bmatrix}U^T\\ U_0^T \end{bmatrix}\begin{bmatrix}U & U_0 \end{bmatrix} \ = \ I_n
\end{equation}

\noindent
which implies that 
\begin{equation}
    U_0^TU = 0_{(n-p)\times p} \ , \ U^TU_0 = 0_{p\times (n-p)}. 
\end{equation}

\noindent
Using $A = UDV^T$ in \eqref{eq:CaseI.1.TraceClass}, and
standard matrix algebra leveraging the various orthogonality properties discussed above, we get 
\begin{align}
    & \left(2W^T\Lambda W + I_p\right) - 4 W^T\Lambda\left(\Lambda + \Lambda  W \left(W^T\Lambda W + I_p\right)^{-1}W^T\Lambda\right)^{-1}\Lambda  W\\
    &= \left(2A^TA + I_p\right) - 4 A^T\left(I_n + A \left(A^TA + I_p\right)^{-1}A^T\right)^{-1}A\\
    %&= \left(2VD^2V^T + I_p\right) - 4 VDU^T\left(I_n + UDV^T \left(VD^2V^T + I_p\right)^{-1}VDU^T\right)^{-1}UDV^T\\
    %&= \left(2VD^2V^T + I_p\right) - 4 VDU^T\left(I_n + I_nUDV^T \left(VD^2V^T + I_p\right)^{-1}VDU^TI_n\right)^{-1}UDV^T\\
    &= \left(2VD^2V^T + I_p\right)\\
    &\quad - 4 VDU^T\Bigg(\begin{bmatrix}U & U_0 \end{bmatrix}\begin{bmatrix}U^T\\ U_0^T \end{bmatrix}\\
    &\qquad\qquad\qquad + \begin{bmatrix}U & U_0 \end{bmatrix}\begin{bmatrix}U^T\\ U_0^T \end{bmatrix}UDV^T \left(VD^2V^T + I_p\right)^{-1}VDU^T\begin{bmatrix}U & U_0 \end{bmatrix}\begin{bmatrix}U^T\\ U_0^T \end{bmatrix}\Bigg)^{-1}UDV^T\\
    %&= \left(2VD^2V^T + I_p\right) - 4 VDU^T\Bigg(\begin{bmatrix}U & U_0 \end{bmatrix}\begin{bmatrix}U^T\\ U_0^T \end{bmatrix}\\
    %&\qquad\qquad\qquad + \begin{bmatrix}U & U_0 \end{bmatrix}\begin{bmatrix}I_p\\ 0 \end{bmatrix}DV^T \left(VD^2V^T + I_p\right)^{-1}VD\begin{bmatrix}I_p & 0 \end{bmatrix}\begin{bmatrix}U^T\\ U_0^T \end{bmatrix}\Bigg)^{-1}UDV^T\\
    &= \left(2VD^2V^T + I_p\right)\\
    &\quad - 4 VDU^T\Bigg(\begin{bmatrix}U & U_0 \end{bmatrix}\begin{bmatrix}U^T\\ U_0^T \end{bmatrix}\\
    &\qquad\qquad\qquad + \begin{bmatrix}U & U_0 \end{bmatrix}\begin{bmatrix}DV^T \left(VD^2V^T + I_p\right)^{-1}VD & 0 \\0 & 0 \end{bmatrix}\begin{bmatrix}U^T\\ U_0^T \end{bmatrix}\Bigg)^{-1}UDV^T\\
    &= \left(2VD^2V^T + I_p\right)\\
    &\quad - 4 VDU^T\begin{bmatrix}U & U_0 \end{bmatrix}\Bigg(I_n + \begin{bmatrix}DV^T \left(VD^2V^T + I_p\right)^{-1}VD & 0 \\0 & 0 \end{bmatrix}\Bigg)^{-1}\begin{bmatrix}U^T\\ U_0^T \end{bmatrix}UDV^T\\
    %&= \left(2VD^2V^T + I_p\right) - 4 VD\begin{bmatrix}I_p & 0 \end{bmatrix}\Bigg(I_n + \begin{bmatrix}DV^T \left(VD^2V^T + I_p\right)^{-1}VD & 0 \\0 & 0 \end{bmatrix}\Bigg)^{-1}\begin{bmatrix}I_p\\ 0 \end{bmatrix}DV^T\\
    %&= \left(2VD^2V^T + I_p\right) - 4 VD\begin{bmatrix}I_p & 0 \end{bmatrix}\begin{bmatrix}\left(I_p + DV^T \left(VD^2V^T + I_p\right)^{-1}VD\right) & 0 \\0 & I_{n-p} \end{bmatrix}^{-1}\begin{bmatrix}I_p\\ 0 \end{bmatrix}DV^T\\
    &= \left(2VD^2V^T + I_p\right) - 4 VD\begin{bmatrix}I_p & 0 \end{bmatrix}\begin{bmatrix}\left(I_p + DV^T \left(VD^2V^T + I_p\right)^{-1}VD\right)^{-1} & 0 \\0 & I_{n-p} \end{bmatrix}\begin{bmatrix}I_p\\ 0 \end{bmatrix}DV^T\\
    &= \left(2VD^2V^T + I_p\right) - 4 VD\left(I_p + DV^T \left(VD^2V^T + I_p\right)^{-1}VD\right)^{-1}DV^T\\
    %&= \left(2VD^2V^T + VV^T\right) - 4 VD\left(I_p + DV^T \left(VD^2V^T + VV^T\right)^{-1}VD\right)^{-1}DV^T\\
    &= \left(2VD^2V^T + VV^T\right) - 4 VD\left(I_p + D \left(D^2 + I_p\right)^{-1}D\right)^{-1}DV^T\\
    &= V\left(2D^2 + I_p\right)V^T - 4 VD\left(I_p + \frac{D^2}{D^2 + I_p}\right)^{-1}DV^T\\
    %&= V\left(2D^2 + I_p\right)V^T - 4 VD\left(\frac{2D^2 + I_p}{D^2 + I_p}\right)^{-1}DV^T\\
    %&= V\left(2D^2 + I_p\right)V^T - 4 VD\left(\frac{D^2 + I_p}{2D^2 + I_p}\right)DV^T\\
    %&= V\left[\left(2D^2 + I_p\right) - 4 D\left(\frac{D^2 + I_p}{2D^2 + I_p}\right)D\right]V^T\\
    %&= V\left[\frac{\left(2D^2 + I_p\right)^2 - 4\left(D^4 + D^2\right)}{2D^2 + I_p}\right]V^T\\
    %&= V\left[\frac{4D^4 + 4D^2 + I_p - 4\left(D^4 + D^2\right)}{2D^2 + I_p}\right]V^T\\
    &= V\left[\frac{I_p}{2D^2 + I_p}\right]V^T
    \label{eq:CaseI.2.TraceClass}
\end{align}

\medskip

\noindent
\textbf{\underline{High-dimensional setting}:} When 
$n < p$, the Singular Value Decomposition of $A = \Lambda^{1/2} W$ can be written as
\begin{equation}
    A_{n\times p} = V_{n\times n}D_{n\times n}U_{n\times p}^T
    \label{eq:SVD.A.CaseII.TraceClass}
\end{equation}
where, $U_{p\times n}$ is a semi-orthogonal matrix i.e. $U^TU=I_n$, $V$ is an orthogonal matrix i.e. $V^TV = VV^T = I_n$, and $D$ is a diagonal matrix with singular values of $A$ being the diagonal entries. Since, $U_{p\times n}$ is a semi-orthogonal matrix with full column rank $n$, there exists a matrix $U_{{0}_{ \ p \times (p-n)}}$ such that the matrix $\begin{bmatrix}U_{p\times n} & U_{{0}_{ \ p \times (p-n)}} \end{bmatrix}$ becomes an orthogonal matrix, i.e.,
\begin{equation}
    \begin{bmatrix}U & U_0 \end{bmatrix}\begin{bmatrix}U^T\\ U_0^T \end{bmatrix} \ = \ \begin{bmatrix}U^T\\ U_0^T \end{bmatrix}\begin{bmatrix}U & U_0 \end{bmatrix} \ = \ I_p
\end{equation}

\noindent
which implies that 
\begin{equation}
    U_0^TU = 0_{(p-n)\times n} \ , \ U^TU_0 = 0_{n\times (p-n)}. 
\end{equation}

\noindent
Again, using $A = VDU^T$ in \eqref{eq:CaseI.1.TraceClass}, 
and standard matrix algebra leveraging the various 
orthogonality properties discussed above, we
get
\begin{align}
    & \left(2W^T\Lambda W + I_p\right) - 4 W^T\Lambda\left(\Lambda + \Lambda  W \left(W^T\Lambda W + I_p\right)^{-1}W^T\Lambda\right)^{-1}\Lambda  W\\
    &= \left(2A^TA + I_p\right) - 4 A^T\left(I_n + A \left(A^TA + I_p\right)^{-1}A^T\right)^{-1}A\\
    %&= \left(2UD^2U^T + I_p\right) - 4 UDV^T\left(I_n + VDU^T \left(UD^2U^T + I_p\right)^{-1}UDV^T\right)^{-1}VDU^T\\
    %&= \left(2UD^2U^T + I_p\right) - 4 UDV^T\left(VV^T + VDU^T \left(UD^2U^T + I_p\right)^{-1}UDV^T\right)^{-1}VDU^T\\
    &= \left(2UD^2U^T + I_p\right) - 4 UD\left(I_n + DU^T \left(UD^2U^T + I_p\right)^{-1}UD\right)^{-1}DU^T\\
    &= \left(2UD^2U^T + I_p\right) - 4 UD\left(I_n + DU^T \left(UD^2U^T + \begin{bmatrix}U & U_0 \end{bmatrix}\begin{bmatrix}U^T\\ U_0^T \end{bmatrix}\right)^{-1}UD\right)^{-1}DU^T\\
    %&= \left(2UD^2U^T + I_p\right) - 4 UD\left(I_n + DU^T \left(UD^2U^T + UU^T + U_0U_0^T\right)^{-1}UD\right)^{-1}DU^T\\
    &= \left(2UD^2U^T + I_p\right) - 4 UD\left(I_n + DU^T \left(U\left(D^2 + I_n\right)U^T + U_0U_0^T\right)^{-1}UD\right)^{-1}DU^T\\
    &= \left(2UD^2U^T + I_p\right) - 4 UD\left(I_n + DU^T \left(U\left(D^2 + I_n\right)^{-1}U^T + U_0U_0^T\right)UD\right)^{-1}DU^T\\
    %&= \left(2UD^2U^T + I_p\right) - 4 UD\left(I_n + D\left(D^2 + I_n\right)^{-1}D\right)^{-1}DU^T\\
    %&= \left(2UD^2U^T + UU^T + U_0U_0^T\right) - 4 UD\left(I_n + D\left(D^2 + I_n\right)^{-1}D\right)^{-1}DU^T\\
    %&= U\left(2D^2 + I_n\right)U^T + U_0U_0^T - U\left[4D^2\left(I_n + D^2\left(D^2+I_n\right)^{-1}\right)^{-1}\right]U^T\\
    &= U\left(2D^2 + I_n\right)U^T + U_0U_0^T - U\left[4D^2\left(D^2 + I_n\right)\left(2D^2+I_n\right)^{-1}\right]U^T\\
    %&= U\left(2D^2 + I_n\right)U^T + U_0U_0^T - U\left[\frac{4D^2\left(D^2+I_n\right)}{2D^2+I_n}\right]U^T\\
    &= U\left[\left(2D^2 + I_n\right) - \frac{4D^2\left(D^2+I_n\right)}{2D^2+I_n}\right]U^T + U_0U_0^T\\
    &= U\left[\frac{I_n}{2D^2+I_n}\right]U^T + U_0U_0^T
    \label{eq:CaseII.2.TraceClass}
\end{align}
Now, if we denote
\begin{equation}
    \Omega\left(\Lambda\right) := \left(2W^T\Lambda W + I_p\right) - 4 W^T\Lambda\left(\Lambda + \Lambda  W \left(W^T\Lambda W + I_p\right)^{-1}W^T\Lambda\right)^{-1}\Lambda  W
\end{equation}
then from \eqref{eq:CaseI.2.TraceClass} and \eqref{eq:CaseII.2.TraceClass}, it follows that
\begin{equation}
    \Omega\left(\Lambda\right) =
\left\{
	\begin{array}{ll}
		V\left[\frac{I_p}{2D^2 + I_p}\right]V^T  & \mbox{if } n \geq p \\
		\\U\left[\frac{I_n}{2D^2+I_n}\right]U^T + U_0U_0^T & \mbox{if } n < p
	\end{array}
\right.
\label{eq:Omega.Lambda.TraceClass}
\end{equation}  
Then, clearly $\Omega\left(\Lambda\right)$ is a positive definite matrix for both the cases $n\geq p$ and $n<p$, and for all $\bl \in\R^n_{+}$. This implies that \ $G\left({\boldsymbol \theta},\bl\right) = {\boldsymbol \theta}^T\Omega\left(\Lambda\right){\boldsymbol \theta}$ \ is a positive definite quadratic form in ${\boldsymbol \theta}$ for all ${\boldsymbol \theta}\in\R^p$ and $\bl\in\R^n_{+}$. Moreover,
\begin{equation}
    \Sigma\left(\Lambda\right) := \Omega\left(\Lambda\right)^{-1} =
\left\{
	\begin{array}{ll}
		V\left(2D^2 + I_p\right)V^T  & \mbox{if } n \geq p \\
		\\U\left(2D^2+I_n\right)U^T + U_0U_0^T & \mbox{if } n < p
	\end{array}
\right.
\label{eq:Sigma.Lambda.TraceClass}
\end{equation}
Now, from \eqref{eq:integral.value.5.TraceClass} we get
\begin{align}
    I \ &\leq \ C_7\int_{\R^n_{+}}\int_{\R^p}\text{max}\left\{1,\lambda_{\text{max}}^{p/2}\right\} \times \left[\prod_{i=1}^{n}\lambda_i^{\frac{\nu}{2} - 1}\right]\times \  \exp\left[-\frac{\nu}{2}\sum_{i=1}^{n}\lambda_i\right]\\
    & \qquad\qquad\qquad\qquad \times \left[C_3 \ + \ C_6\left(\left({\boldsymbol \theta}+\tilde{\bf c}\right)^T\left({\boldsymbol \theta}+\tilde{\bf c}\right)\right)^{\frac{n\nu}{2}}\right]\\
    & \qquad\qquad\qquad\qquad \times \exp\left[-\frac{1}{2}\Big\{G\left({\boldsymbol \theta}, \bl\right)\Big\}\right] \ d{\boldsymbol \theta} \ d\bl\\
    &= \ C_7\int_{\R^n_{+}}\int_{\R^p}\text{max}\left\{1,\lambda_{\text{max}}^{p/2}\right\} \times \left[\prod_{i=1}^{n}\lambda_i^{\frac{\nu}{2} - 1}\right]\times \  \exp\left[-\frac{\nu}{2}\sum_{i=1}^{n}\lambda_i\right]\\
    & \qquad\qquad\qquad\qquad \times \left[C_3 \ + \ C_6\left(\left({\boldsymbol \theta}+\tilde{\bf c}\right)^T\left({\boldsymbol \theta}+\tilde{\bf c}\right)\right)^{\frac{n\nu}{2}}\right]\\
    & \qquad\qquad\qquad\qquad \times \exp\left[-\frac{1}{2} \ {\boldsymbol \theta}^T\Omega\left(\Lambda\right){\boldsymbol \theta}\right] \ d{\boldsymbol \theta} \ d\bl\\
    &= \ C_7\int_{\R^n_{+}}\text{max}\left\{1,\lambda_{\text{max}}^{p/2}\right\} \times \left[\prod_{i=1}^{n}\lambda_i^{\frac{\nu}{2} - 1}\right]\times \  \exp\left[-\frac{\nu}{2}\sum_{i=1}^{n}\lambda_i\right]\\
    & \qquad\qquad \times \left(\int_{\R^p}\left[C_3 \ + \ C_6\left(\left({\boldsymbol \theta}+\tilde{\bf c}\right)^T\left({\boldsymbol \theta}+\tilde{\bf c}\right)\right)^{\frac{n\nu}{2}}\right] \exp\left[-\frac{1}{2} \ {\boldsymbol \theta}^T\Omega\left(\Lambda\right){\boldsymbol \theta}\right] \ d{\boldsymbol \theta}\right) \ d\bl\\
    & \label{eq:integral.value.6.TraceClass}
\end{align}

\noindent
{\bf Step VI: An upper bound for the inner integral in 
\eqref{eq:integral.value.6.TraceClass}.} We now derive 
an upper bound for the inner integral in 
\eqref{eq:integral.value.6.TraceClass} using properties
of the multivariate normal distribution. Note that 
\begin{align}
    & \ \int_{\R^p}\left[C_3 \ + \ C_6\left(\left({\boldsymbol \theta}+\tilde{\bf c}\right)^T\left({\boldsymbol \theta}+\tilde{\bf c}\right)\right)^{\frac{n\nu}{2}}\right] \exp\left[-\frac{1}{2} \ {\boldsymbol \theta}^T\Omega\left(\Lambda\right){\boldsymbol \theta}\right] \ d{\boldsymbol \theta}\\
    &= \ \left(2\pi\right)^{p/2}\sqrt{\text{det}\left(\Sigma\left(\Lambda\right)\right)}\int_{\R^p}\left(2\pi\right)^{-p/2} \ \text{det}\left(\Sigma\left(\Lambda\right)\right)^{-1/2}\left[C_3 \ + \ C_6\left(\left({\boldsymbol \theta}+\tilde{\bf c}\right)^T\left({\boldsymbol \theta}+\tilde{\bf c}\right)\right)^{\frac{n\nu}{2}}\right]\\
    &\qquad\qquad\qquad\qquad\qquad\qquad\qquad\qquad\qquad\qquad\qquad \times \exp\left[-\frac{1}{2} \ {\boldsymbol \theta}^T\Omega\left(\Lambda\right){\boldsymbol \theta}\right] \ d{\boldsymbol \theta}\\
    &= \ \left(2\pi\right)^{p/2}\sqrt{\text{det}\left(\Sigma\left(\Lambda\right)\right)} \ \times \ E_{\mathcal{N}_p\left({\bm 0}, \ \Sigma\left(\Lambda\right)\right)}\left(\left[C_3 \ + \ C_6\left(\left({\boldsymbol \theta}+\tilde{\bf c}\right)^T\left({\boldsymbol \theta}+\tilde{\bf c}\right)\right)^{\frac{n\nu}{2}}\right]\right)
    \label{eq:integral.value.6.inner.integral.TraceClass}
\end{align}
where, $\mathcal{N}_p\left({\bm 0}, \ \Sigma\left(\Lambda\right)\right)$ stands for multivariate normal distribution with mean vector ${\bm 0}$ and covariance matrix $\Sigma\left(\Lambda\right)$. Observe that 
\begin{align}
    & \ E_{\mathcal{N}_p\left({\bm 0}, \ \Sigma\left(\Lambda\right)\right)}\left(\left[C_3 \ + \ C_6\left(\left({\boldsymbol \theta}+\tilde{\bf c}\right)^T\left({\boldsymbol \theta}+\tilde{\bf c}\right)\right)^{\frac{n\nu}{2}}\right]\right)\\
    %&= \ C_3 \ + \ C_6 \ E_{\text{MVN}\left(0,\Sigma\left(\Lambda\right)\right)}\left[\left(\left({\boldsymbol \theta}+\tilde{c}\right)^T\left({\boldsymbol \theta}+\tilde{c}\right)\right)^{\frac{n\nu}{2}}\right]\\
    &= \ C_3 \ + \ C_6 \ E_{\mathcal{N}_p\left({\bm 0}, \ \Sigma\left(\Lambda\right)\right)}\left[\left(\left({\boldsymbol \theta}+\tilde{\bf c}\right)^T\Omega\left(\Lambda\right)^{1/2} \Sigma\left(\Lambda\right)\Omega\left(\Lambda\right)^{1/2} \left({\boldsymbol \theta}+\tilde{\bf c}\right)\right)^{\frac{n\nu}{2}}\right]\\
    &\leq \ C_3 \ + \ C_6 \ \text{eig}_{\text{max}}^{n\nu/2}\left(\Sigma\left(\Lambda\right)\right) \  E_{\mathcal{N}_p\left({\bm 0}, \ \Sigma\left(\Lambda\right)\right)}\left[\left(\left({\boldsymbol \theta}+\tilde{\bf c}\right)^T\Omega\left(\Lambda\right) \left({\boldsymbol \theta}+\tilde{\bf c}\right)\right)^{\frac{n\nu}{2}}\right]\\
    &= \ C_3 \ + \ C_6 \ \text{eig}_{\text{max}}^{n\nu/2}\left(\Sigma\left(\Lambda\right)\right) \  E_{\mathcal{N}_p\left({\bm 0}, \ \Sigma\left(\Lambda\right)\right)}\left[\left({\boldsymbol \theta}^T\Omega\left(\Lambda\right){\boldsymbol \theta} + 2\tilde{\bf c}^T\Omega\left(\Lambda\right){\boldsymbol \theta} + \tilde{\bf c}^T\Omega\left(\Lambda\right)\tilde{\bf c}\right)^{\frac{n\nu}{2}}\right]\\
    & \label{eq:integral.value.6_1.inner.integral.TraceClass}
\end{align}

\noindent
since $\Sigma\left(\Lambda\right) \ \preceq \ 
\text{eig}_{\text{max}}\left(\Sigma\left(\Lambda\right)
\right) I_p$, and 
\begin{align}
    & \ E_{\mathcal{N}_p\left({\bm 0}, \ \Sigma\left(\Lambda\right)\right)}\left[\left({\boldsymbol \theta}^T\Omega\left(\Lambda\right){\boldsymbol \theta} + 2\tilde{\bf c}^T\Omega\left(\Lambda\right){\boldsymbol \theta} + \tilde{\bf c}^T\Omega\left(\Lambda\right)\tilde{\bf c}\right)^{\frac{n\nu}{2}}\right]\\
    &\leq \ E_{\mathcal{N}_p\left({\bm 0}, \ \Sigma\left(\Lambda\right)\right)}\left[\left({\boldsymbol \theta}^T\Omega\left(\Lambda\right){\boldsymbol \theta} + 2\left|\tilde{\bf c}^T\Omega\left(\Lambda\right){\boldsymbol \theta}\right| + \tilde{\bf c}^T\Omega\left(\Lambda\right)\tilde{\bf c}\right)^{\frac{n\nu}{2}}\right]\\
    &\leq \ 3^{n\nu/2} \ E_{\mathcal{N}_p\left({\bm 0}, \ \Sigma\left(\Lambda\right)\right)}\left[\left({\boldsymbol \theta}^T\Omega\left(\Lambda\right){\boldsymbol \theta}\right)^{\frac{n\nu}{2}} + \left(2\left|\tilde{\bf c}^T\Omega\left(\Lambda\right){\boldsymbol \theta}\right|\right)^{\frac{n\nu}{2}} + \left(\tilde{\bf c}^T\Omega\left(\Lambda\right)\tilde{\bf c}\right)^{\frac{n\nu}{2}}\right]\\
    & \ \Big[\text{since, for any }a,b,c,n \in \R_{+}\cup \{0\}, \ \left(a+b+c\right)^n \leq 3^n\left(a^n+b^n+c^n\right)\Big]\\
    &= \ 3^{n\nu/2} \ \Bigg[E_{\mathcal{N}_p\left({\bm 0}, \ \Sigma\left(\Lambda\right)\right)}\left[\left({\boldsymbol \theta}^T\Omega\left(\Lambda\right){\boldsymbol \theta}\right)^{\frac{n\nu}{2}}\right] + E_{\mathcal{N}_p\left({\bm 0}, \ \Sigma\left(\Lambda\right)\right)}\left[\left(2\left|\tilde{\bf c}^T\Omega\left(\Lambda\right){\boldsymbol \theta}\right|\right)^{\frac{n\nu}{2}}\right]\\
    &\qquad\qquad\qquad\qquad\qquad\qquad + E_{\mathcal{N}_p\left({\bm 0}, \ \Sigma\left(\Lambda\right)\right)}\left[\left(\tilde{\bf c}^T\Omega\left(\Lambda\right)\tilde{\bf c}\right)^{\frac{n\nu}{2}}\right]\Bigg]
    \label{eq:integral.value.6_2.inner.integral.TraceClass}
\end{align}

\noindent
We will show that each term in \eqref{eq:integral.value.6_2.inner.integral.TraceClass} is uniformly bounded in $\bl$. Note that if ${\boldsymbol \theta} \ \sim \ \mathcal{N}_p\left({\bm 0}, \ \Sigma\left(\Lambda\right)\right)$, then ${\boldsymbol\theta}^T\Omega\left(\Lambda\right){\boldsymbol \theta} \ \sim \ \chi^2_p$. Hence 
$E_{\mathcal{N}_p\left({\bm 0}, \ \Sigma\left(\Lambda\right)\right)}\left[\left({\boldsymbol \theta}^T\Omega\left(\Lambda\right){\boldsymbol \theta}\right)^{\frac{n\nu}{2}}\right]$ is a finite quantity free of $\bl$. 

Again, using the fact that $2\left|{\bf a}^T{\bf b}\right|\leq {\bf a}^T{\bf a} + {\bf b}^T{\bf b}$, and taking ${\bf a} = \Omega\left(\Lambda\right)^{\frac{1}{2}}\tilde{\bf c}$ and ${\bf b} = \Omega\left(\Lambda\right)^{\frac{1}{2}}{\boldsymbol \theta}$, we get
\begin{align}
    2\left|\tilde{\bf c}^T\Omega\left(\Lambda\right){\boldsymbol \theta}\right| \ &\leq \ \tilde{\bf c}^T\Omega\left(\Lambda\right)\tilde{\bf c} \ + \ {\boldsymbol \theta}^T\Omega\left(\Lambda\right){\boldsymbol \theta}\\
    &\leq \ \text{eig}_{\text{max}}\left(\Omega\left(\Lambda\right)\right) \ \tilde{\bf c}^T\tilde{\bf c} \ + \ {\boldsymbol \theta}^T\Omega\left(\Lambda\right){\boldsymbol \theta}\\
    & \ \Big[\text{since, }\Omega\left(\Lambda\right) \ \preceq \ \text{eig}_{\text{max}}\left(\Omega\left(\Lambda\right)\right) I_p\Big]
\end{align}
However, from the expression of $\Omega\left(\Lambda\right)$ in \eqref{eq:Omega.Lambda.TraceClass}, it is easy to see that $\Omega\left(\Lambda\right) \preceq I_p$ and  hence $\text{eig}_{\text{max}}\left(\Omega\left(\Lambda\right)\right) \leq 1$. It follows that 
\begin{align}
     E_{\mathcal{N}_p\left({\bm 0}, \ \Sigma\left(\Lambda\right)\right)}\left[\left(2\left|\tilde{\bf c}^T\Omega\left(\Lambda\right){\boldsymbol \theta}\right|\right)^{\frac{n\nu}{2}}\right] \ &\leq \ 2^{n\nu/2} \ E_{\mathcal{N}_p\left({\bm 0}, \ \Sigma\left(\Lambda\right)\right)}\left[\left(\tilde{\bf c}^T\tilde{\bf c}\right)^{\frac{n\nu}{2}} + \left({\boldsymbol \theta}^T\Omega\left(\Lambda\right){\boldsymbol \theta}\right)^{\frac{n\nu}{2}}\right]\\
    &= \ 2^{n\nu/2} \ \left[\left(\tilde{\bf c}^T\tilde{\bf c}\right)^{\frac{n\nu}{2}} + E_{\mathcal{N}_p\left({\bm 0}, \ \Sigma\left(\Lambda\right)\right)}\left[\left({\boldsymbol \theta}^T\Omega\left(\Lambda\right){\boldsymbol \theta}\right)^{\frac{n\nu}{2}}\right]\right]\\
    %&= \ 2^{n\nu/2} \ \left[\left(\tilde{c}^T\tilde{c}\right)^{\frac{n\nu}{2}} + C_8\right] \qquad \left(\text{by }\eqref{eq:integral.value.6_2.1.inner.integral.TraceClass}\right)\\
    &= \ C_9 \quad \left(\text{say}\right)
    \label{eq:integral.value.6_2.2.inner.integral.TraceClass}
\end{align}

\noindent
where $C_9$ is finite and independent of $\bl$ based on
the observations above. Finally, since all eigenvalues 
of $\Omega\left(\Lambda\right)$ are non-negative and  
bounded above by $1$, it follows that 
$$
\left(\tilde{\bf c}^T\Omega\left(\Lambda\right)\tilde{\bf c}\right)^{\frac{n\nu}{2}} \leq \left(\tilde{\bf c}^T \tilde{\bf c}\right)^{\frac{n\nu}{2}}. 
$$

\noindent
Using \eqref{eq:integral.value.6_1.inner.integral.TraceClass} and \eqref{eq:integral.value.6_2.inner.integral.TraceClass}, we get
\begin{align}
    & \ E_{\mathcal{N}_p\left({\bm 0}, \ \Sigma\left(\Lambda\right)\right)}\left(\left[C_3 \ + \ C_6\left(\left({\boldsymbol \theta}+\tilde{\bf c}\right)^T\left({\boldsymbol \theta}+\tilde{\bf c}\right)\right)^{\frac{n\nu}{2}}\right]\right)\\
    %&\leq \ C_3 \ + \ C_6 \ \text{eig}_{\text{max}}^{n\nu/2}\left(\Sigma\left(\Lambda\right)\right) \  E_{\text{MVN}\left(0,\Sigma\left(\Lambda\right)\right)}\left[\left({\boldsymbol \theta}^T\Omega\left(\Lambda\right){\boldsymbol \theta} + 2\tilde{c}^T\Omega\left(\Lambda\right){\boldsymbol \theta} + \tilde{c}^T\Omega\left(\Lambda\right)\tilde{c}\right)^{\frac{n\nu}{2}}\right]\\
    &\leq \ C_3 \ + \ C_{12} \ \text{eig}_{\text{max}}^{n\nu/2}\left(\Sigma\left(\Lambda\right)\right)\\
    &= \ C_3 \ + \ C_{12}\left(2 d_{\text{max}}^2 + 1\right)^{\frac{n\nu}{2}} \qquad \left(\text{by the expression of }\Sigma\left(\Lambda\right) \ \text{in }\eqref{eq:Sigma.Lambda.TraceClass}\right)
    \label{eq:integral.value.6_4.inner.integral.TraceClass}
\end{align}

\noindent
where $C_{12}$ is an appropriate constant independent 
of $\bl$ and $d_{\text{max}}$ is the largest element of the diagonal matrix $D$ in the expression of \eqref{eq:Sigma.Lambda.TraceClass}. Plugging this in  \eqref{eq:integral.value.6.inner.integral.TraceClass}, we get the following upper bound for the inner integral in \eqref{eq:integral.value.6.TraceClass}. 
\begin{align}
    & \ \int_{\R^p}\left[C_3 \ + \ C_6\left(\left({\boldsymbol \theta}+\tilde{\bf c}\right)^T\left({\boldsymbol \theta}+\tilde{\bf c}\right)\right)^{\frac{n\nu}{2}}\right] \exp\left[-\frac{1}{2} \ {\boldsymbol \theta}^T\Omega\left(\Lambda\right){\boldsymbol \theta}\right] \ d{\boldsymbol \theta}\\
    %&= \ \left(2\pi\right)^{p/2}\sqrt{\text{det}\left(\Sigma\left(\Lambda\right)\right)} \ \times \ E_{\text{MVN}\left(0,\Sigma\left(\Lambda\right)\right)}\left(\left[C_3 \ + \ C_6\left(\left({\boldsymbol \theta}+\tilde{c}\right)^T\left({\boldsymbol \theta}+\tilde{c}\right)\right)^{\frac{n\nu}{2}}\right]\right)\\
    &\leq \ \left(2\pi\right)^{p/2}\sqrt{\text{det}\left(\Sigma\left(\Lambda\right)\right)} \ \left[C_3 \ + \ C_{12}\left(2 d_{\text{max}}^2 + 1\right)^{\frac{n\nu}{2}}\right]. 
    \label{eq:integral.value.6.inner.integral.upper.bound.TraceClass}
\end{align}

\noindent
This in turn leads to the following upper bound for the
integral $I$ of interest. 
\begin{align}
    I \ &\leq \ C_7\left(2\pi\right)^{p/2}\int_{\R^n_{+}}\text{max}\left\{1,\lambda_{\text{max}}^{p/2}\right\} \times \left[\prod_{i=1}^{n}\lambda_i^{\frac{\nu}{2} - 1}\right]\times \  \exp\left[-\frac{\nu}{2}\sum_{i=1}^{n}\lambda_i\right]\\
    & \qquad\qquad\qquad\qquad \times \sqrt{\text{det}\left(\Sigma\left(\Lambda\right)\right)} \ \left[C_3 \ + \ C_{12}\left(2 d_{\text{max}}^2 + 1\right)^{\frac{n\nu}{2}}\right] \ d\bl. 
    \label{eq:integral.value.7.TraceClass}
\end{align}

\noindent
Our final steps will be to bound the (single) integral 
in \eqref{eq:integral.value.7.TraceClass} separately in
the low and high-dimensional settings. 

\medskip
\noindent
{\bf Step VII: An upper bound for the integral in 
\eqref{eq:integral.value.7.TraceClass} in the low-dimensional setting.} 
Suppose $n \geq p$. By the expression of $\Sigma\left(\Lambda\right)$ in 
\eqref{eq:Sigma.Lambda.TraceClass} for $n\geq p$, we have
\begin{align}
    \text{det}\left(\Sigma\left(\Lambda\right)\right) \ = \ \text{det}\left(V\left(2D^2 + I_p\right)V^T\right)
    = \ \prod_{i=1}^{p}\left(2d_i^2+1\right)
    \leq \ \left(2d_{\text{max}}^2+1\right)^p, \\
    \label{eq:CaseI.determinant.Sigma.Lambda.TraceClass}
\end{align}

\noindent
where $\{d_i\}_{i=1}^p$ are the singular values of $A$,
and $d_{\text{max}}$ is the largest singular value of 
$A = \Lambda^{1/2}W$ in \eqref{eq:SVD.A.CaseI.TraceClass}. Note that 
\begin{align}
    d_{\text{max}}^2 \ = \ \text{eig}_{\text{max}}\left(AA^T\right) 
    &= \ \text{eig}_{\text{max}}\left(\Lambda^{1/2}WW^T\Lambda^{1/2}\right). 
    \label{eq:CaseI.d_max.squared.1.TraceClass}
\end{align}
Since $WW^T$ is a fixed positive semi-definite matrix, there exists a large positive real number $C_{13}$ such that
\begin{align}
\Lambda^{1/2}WW^T\Lambda^{1/2} \preceq C_{13}\Lambda 
\preceq C_{13}\lambda_{\text{max}} I_n, 
\end{align}
where $\lambda_{\text{max}}$ is the largest element in the matrix $\Lambda$. It follows from \eqref{eq:CaseI.d_max.squared.1.TraceClass} that $
d_{\text{max}}^2 \leq C_{13} \lambda_{\text{max}}$.  Using \eqref{eq:integral.value.7.TraceClass}, 
\eqref{eq:CaseI.determinant.Sigma.Lambda.TraceClass} 
and the $c_r$-inequality, we get
\begin{align}
    I \ &\leq \ C_7\left(2\pi\right)^{p/2}\int_{\R^n_{+}}\text{max}\left\{1,\lambda_{\text{max}}^{p/2}\right\} \times \left[\prod_{i=1}^{n}\lambda_i^{\frac{\nu}{2} - 1}\right]\times \  \exp\left[-\frac{\nu}{2}\sum_{i=1}^{n}\lambda_i\right]\\
    & \qquad\qquad\qquad\qquad \times \left(2d_{\text{max}}^2+1\right)^{\frac{p}{2}} \ \left[C_3 \ + \ C_{12}\left(2 d_{\text{max}}^2 + 1\right)^{\frac{n\nu}{2}}\right] \ d\bl\\
    &\leq \ C_7\left(2\pi\right)^{p/2}\int_{\R^n_{+}}\text{max}\left\{1,\lambda_{\text{max}}^{p/2}\right\} \times \left[\prod_{i=1}^{n}\lambda_i^{\frac{\nu}{2} - 1}\right]\times \  \exp\left[-\frac{\nu}{2}\sum_{i=1}^{n}\lambda_i\right]\\
    & \qquad\qquad\qquad\qquad \times \left(2C_{13}\lambda_{\text{max}}+1\right)^{\frac{p}{2}} \ \left[C_3 \ + \ C_{12}\left(2C_{13}\lambda_{\text{max}}+1\right)^{\frac{n\nu}{2}}\right] \ d\bl\\
    &\leq \ C_7\left(2\pi\right)^{p/2}\int_{\R^n_{+}}\left(1 +\lambda_{\text{max}}^{p/2}\right) \times \left[\prod_{i=1}^{n}\lambda_i^{\frac{\nu}{2} - 1}\right]\times \  \exp\left[-\frac{\nu}{2}\sum_{i=1}^{n}\lambda_i\right]\\
    & \qquad\qquad\qquad\qquad \times 2^{p/2}\left(C_{14}\lambda_{\text{max}}^{p/2}+1\right) \ \left[C_3 \ + \ C_{12}2^{n\nu/2}\left(C_{14}^{'}\lambda_{\text{max}}^{n\nu/2}+1\right)\right] \ d\bl, 
    \end{align}

\noindent
where $C_{14} = \left(2C_{13}\right)^{p/2}$ and  
$C_{14}^{'} = \left(2C_{13}\right)^{n\nu/2}$. 
Expanding the product of the polynomial terms in 
$\lambda_{\text{max}}$ in the integrand gives 
\begin{align}
I &\leq \ C_{15}\int_{\R^n_{+}}\left(\lambda_{\text{max}}^{p + \frac{n\nu}{2}} + \lambda_{\text{max}}^p + \lambda_{\text{max}}^{\left(p+n\nu\right)/2} + \lambda_{\text{max}}^{p/2} + \lambda_{\text{max}}^{n\nu/2}+1\right)\\
& \qquad\qquad\qquad\qquad \times \left[\prod_{i=1}^{n}\lambda_i^{\frac{\nu}{2} - 1}\right]\times \  \exp\left[-\frac{\nu}{2}\sum_{i=1}^{n}\lambda_i\right] \ d\bl 
\end{align}

\noindent
for an appropriate finite constant $C_{15}$ which does 
not depend on $\bl$. Using the fact that 
$\lambda_{\text{max}}^r \leq \sum_{j=1}^n \lambda_j^r$ 
for any positive $r$, we get 
\begin{align}
I &\leq \ C_{15} \sum_{j=1}^n \int_{\R^n_{+}}\left(\lambda_j^{p + \frac{n\nu}{2}} + \lambda_j^p + \lambda_j^{\left(p+n\nu\right)/2} + \lambda_j^{p/2} + \lambda_j^{n\nu/2}+1\right)\\
& \qquad\qquad\qquad\qquad \times \left[\prod_{i=1}^{n}\lambda_i^{\frac{\nu}{2} - 1}\right]\times \  \exp\left[-\frac{\nu}{2}\sum_{i=1}^{n}\lambda_i\right] \ d\bl\\
&\leq C_{15} \sum_{j=1}^n \left( \int_{\mathbb{R}_+} \left(\lambda_j^{p + \frac{n\nu}{2}} + \lambda_j^p + \lambda_j^{\left(p+n\nu\right)/2} + \lambda_j^{p/2} + \lambda_j^{n\nu/2}+1\right) \lambda_j^{\frac{\nu}{2} - 1} \ \exp\left[-\frac{\nu}{2}\lambda_j\right] 
d \lambda_j \right)\\
& \qquad\qquad\qquad \times \left( \prod_{i=1,i \neq j}^{n} \int_{\mathbb{R}_+} \lambda_i^{\frac{\nu}{2} - 1} \  \exp\left[-\frac{\nu}{2}\lambda_i\right] \ d \lambda_i \right). \label{eq:integral.value.8.CaseI.TraceClass}
\end{align}
    
\noindent
Since all of the terms in the above bound are integrals
over unnormalized gamma densities (with strictly 
positive and finite shape and rate parameters), it 
follows that $I < \infty$. 

\medskip

\noindent
{\bf Step VIII: An upper bound for the integral in \eqref{eq:integral.value.7.TraceClass} in the high-dimensional setting.} Suppose $n < p$. 
By the expression of $\Sigma\left(\Lambda\right)$ in \eqref{eq:Sigma.Lambda.TraceClass} for $n < p$, we have
\begin{align}
    \text{det}\left(\Sigma\left(\Lambda\right)\right) \ &= \ \text{det}\left(U\left(2D^2 + I_n\right)U^T + U_0U_0^T\right)\\
    &= \ \text{det}\left(\begin{bmatrix}
    \left(2D^2 + I_n\right) & 0\\ 0 & I_{p-n}
    \end{bmatrix}\right)\\
    &= \ \prod_{i=1}^{n}\left(2d_i^2+1\right)\\
    &\leq \ \left(2d_{\text{max}}^2+1\right)^n. 
    \label{eq:CaseII.determinant.Sigma.Lambda.TraceClass}
\end{align}

\noindent
By a similar argument as for the $n\geq p$ setting, we 
get $d_{\text{max}}^2 \leq C_{13}'\lambda_{\text{max}}$
for an appropriate constant $C_{13}'$ not depending on 
$\bl$. Using \eqref{eq:integral.value.7.TraceClass}, 
\eqref{eq:CaseII.determinant.Sigma.Lambda.TraceClass} and the $c_r$-inequality, we get
\begin{align}
    I \ &\leq \ C_7\left(2\pi\right)^{p/2}\int_{\R^n_{+}}\text{max}\left\{1,\lambda_{\text{max}}^{p/2}\right\} \times \left[\prod_{i=1}^{n}\lambda_i^{\frac{\nu}{2} - 1}\right]\times \  \exp\left[-\frac{\nu}{2}\sum_{i=1}^{n}\lambda_i\right]\\
    & \qquad\qquad\qquad\qquad \times \left(2d_{\text{max}}^2+1\right)^{\frac{n}{2}} \ \left[C_3 \ + \ C_{12}\left(2 d_{\text{max}}^2 + 1\right)^{\frac{n\nu}{2}}\right] \ d\bl\\
    &\leq \ C_7\left(2\pi\right)^{p/2}\int_{\R^n_{+}}\text{max}\left\{1,\lambda_{\text{max}}^{p/2}\right\} \times \left[\prod_{i=1}^{n}\lambda_i^{\frac{\nu}{2} - 1}\right]\times \  \exp\left[-\frac{\nu}{2}\sum_{i=1}^{n}\lambda_i\right]\\
    & \qquad\qquad\qquad\qquad \times \left(2C_{13}'\lambda_{\text{max}}+1\right)^{\frac{n}{2}} \ \left[C_3 \ + \ C_{12}\left(2C_{13}'\lambda_{\text{max}}+1\right)^{\frac{n\nu}{2}}\right] \ d\bl\\
    &\leq \ C_7\left(2\pi\right)^{p/2}\int_{\R^n_{+}}\left(1 +\lambda_{\text{max}}^{p/2}\right) \times \left[\prod_{i=1}^{n}\lambda_i^{\frac{\nu}{2} - 1}\right]\times \  \exp\left[-\frac{\nu}{2}\sum_{i=1}^{n}\lambda_i\right]\\
    & \qquad\qquad\qquad\qquad \times 2^{n/2}\left(C_{13}^*\lambda_{\text{max}}^{n/2}+1\right) \ \left[C_3 \ + \ C_{12}2^{n\nu/2}\left(C_{14}^{*}\lambda_{\text{max}}^{n\nu/2}+1\right)\right] \ d\bl,  \label{eq:integral.value.8.CaseII.TraceClass}
    \end{align}
    
\noindent
where $C_{13}^* = \left(2C_{13}'\right)^{n/2}$ and 
$ C_{14}^{*} =  \left(2C_{13}'\right)^{n\nu/2}$ similar
to the low-dimensional setting. Expanding the product 
of the polynomial terms in $\lambda_{\text{max}}$ in 
the integrand gives
\begin{align}
    I &\leq \ C_{15}' \int_{\mathbb{R}^n_+} \Bigg( \lambda_{\text{max}}^{\frac{p + n + n\nu}{2}} + \lambda_{\text{max}}^{\frac{n+n\nu}{2}} + \lambda_{\text{max}}^{\frac{p+n}{2}} +
    \lambda_{\text{max}}^{\frac{n}{2}}\\
    &\qquad\qquad\qquad\qquad\qquad +
    \lambda_{\text{max}}^{\left(p+n\nu\right)/2} +  \lambda_{\text{max}}^{n\nu/2}+ \lambda_{\text{max}}^{p/2} + 1\Bigg)\\
    & \qquad\qquad\qquad\qquad \times \left[\prod_{i=1}^{n}\lambda_i^{\frac{\nu}{2} - 1}\right]\times \  \exp\left[-\frac{\nu}{2}\sum_{i=1}^{n}\lambda_i\right] \ d\bl, 
    \end{align}

\noindent
for an appropriate finite constant $C_{15}'$ which does
not depend on $\bl$. Using the fact that 
$\lambda_{\text{max}}^r \leq \sum_{j=1}^n \lambda_j^r$ 
for any positive $r$, and similar arguments regarding 
gamma integrals as in 
\eqref{eq:integral.value.8.CaseI.TraceClass} for the 
low-dimensional $n \geq p$ setting, it follows that $I 
< \infty$ in the high-dimensional $n < p$ setting as 
well. This establishes the trace-class property of the 
DA Markov chain. 
\end{proof}

\noindent
As discussed in the introduction, the trace-class property established 
above implies compactness of the Markov operator $K$, which implies 
that the corresponding DA Markov chain is geometrically ergodic. 
\begin{corollary} \label{cor1}
For $\nu >2$, the DA Markov chain with transition density $k$ (in 
\eqref{eq:markov.transition.density}) is geometrically ergodic for an 
arbitrary choice of the design matrix $X$, sample size $n$, number of 
predictors $p$, prior mean vector $\bb_a$, and (positive
definite) prior precision matrix $\Sigma_a$. 
\end{corollary}

\section{Numerical Illustrations} \label{sec:numerical:illustration}

\noindent
This section presents numerical illustrations to compare/contrast the convergence properties of the robit DA and some other relevant Markov chains. To examine both the low and high-dimensional settings, we consider two real data sets, viz., the Lupus data ($n > p$) from \cite{vanDyk:Meng:2001} and the prostate cancer data ($n < p$) from \cite{chunSparsePartialLeast2010}.  We note at the outset that with a view to the main goal/contribution of this paper, viz., theoretical convergence analysis for (robit) Markov chains, our focus in this section centers entirely around exemplification of the said convergences. Prior elicitation for statistical inference are beyond the scope of this section and the paper; we consider the data-driven Zellner's $g$-prior following \cite{Chakraborty:Khare:2017} in the first example, and independent standard normal priors in the second example to run the respective Markov chains. 

As noted in the Introduction, a trace-class property ensures guaranteed improvements in the convergence of a DA algorithm by \textit{sandwiching} 
(see Corollary \ref{cor2} below). To exemplify/visualize this improvement in the current context we alongside consider a sandwich algorithm obtained by inserting an inexpensive random generation step in between the two steps of Algorithm~\ref{algo:robit_da}. More specifically, we consider the sandwich algorithm from \cite{Roy:2012:robit} which inserts a univariate random gamma generation as presented in  Algorithm~\ref{algo:robit_sandwich}. Note that both the original DA algorithm and its sandwiched version share the same target stationary distribution for $\bb$.     

\begin{algorithm}
\caption{$(m+1)$-st Iteration of the Robit Sandwich Algorithm}
\label{algo:robit_sandwich}
\begin{enumerate}
    \item Make independent draws from the $Tt_\nu ({\bf x}_i^T \bb^{(m)}, 
    y_i)$ distributions for $1 \leq i \leq n$. Denote the respective draws 
    by $z_1, z_2, \cdots, z_n$. Draw $\lambda_i$ from the Gamma$\left( 
    \frac{\nu+1}{2}, \frac{\nu + \left(z_i - {\bf x}_i^T \bb^{(m)}\right)^2}{2} \right)$ 
    distribution. 
    
    \item Generate $h^2 \sim \operatorname{Gamma}\left(\frac{n}{2}, \frac{{\bf z}^T\Lambda^{1/2}\left(I - Q\right)\Lambda^{1/2}{\bf z}}{2}\right)$, where $Q = \Lambda^{1/2}X \left(X^T \Lambda X + \Sigma_a\right)^{-1}X^T\Lambda^{1/2}$, and subsequently define $z'_i = h z_i$; $1 \leq i \leq n$. 
    
    \item Draw $\bb^{(m+1)}$ from the $\mathcal{N}_p \left( (X^T \Lambda X + \Sigma_a)^{-1} (X^T \Lambda {\bf z}' + \Sigma_a \bb_a) , \ (X^T \Lambda X + \Sigma_a)^{-1} \right)$ distribution. 
\end{enumerate}
\end{algorithm}

\noindent
The trace class property of the robit DA chain 
(Theorem \ref{thm1}) along with results in 
\cite{Khare:Hobert:2011} imply that the following properties 
hold for the sandwich chain. 
\begin{corollary} \label{cor2}
For $\nu >2$, the sandwich Markov chain described in 
Algorithm \ref{algo:robit_sandwich} is trace class (and hence 
geometrically ergodic) for an arbitrary choice of the design 
matrix $X$, sample size $n$, number of predictors $p$, prior 
mean vector $\bb_a$, and (positive definite) prior precision 
matrix $\Sigma_a$. Furthermore, if $\left( \lambda_i 
\right)_{i=0}^\infty$ and $\left( \lambda_i^* 
\right)_{i=0}^\infty$ denote the non-increasing sequences of 
eigenvalues corresponding to the robit DA and sandwich 
operators respectively, then $\lambda_i^* \leq \lambda_i$ for 
every $i \geq 0$, with at least one strict inequality. 
\end{corollary}

\noindent
To facilitate model comparison in each example below we consider two robit models: one with a small degrees of freedom ($\nu = 3$) and one with a large degrees of freedom ($\nu = 1000$). The latter is expected to perform similarly to the probit model, owing to the relationship between a $t$ distribution with large degrees of freedom and the standard normal distribution, which ensures the corresponding likelihoods to be close. To document this similarity herein we also consider the DA and sandwich Markov chains for the probit model from \cite{Chakraborty:Khare:2017}, and note their comparative convergence behaviors vis-a-vis the robit model chains. Collectively we thus consider three models: (a) probit, (b) robit with small $\nu$ and (c) robit with large $\nu$, and for each model we consider two Markov chains -- the original DA chain and a corresponding sandwich chain.

\subsection{Low Dimensional ($n > p$) Setting: Lupus Data Set}

The Lupus data set of \cite{vanDyk:Meng:2001} comprises observations on two antibody molecule predictor variables and a binary outcome variable cataloging occurrences of latent membranous lupus nephritis among $n=55$ patients. Interest lies in regressing the binary outcome on the predictors;  in this regression we include an intercept term which effectively makes the number of predictors to be $p=3$.  The data set has been previously considered in the context of convergence analyses of the probit model DA and sandwich Markov chains \cite{vanDyk:Meng:2001, Roy:Hobert:2007, Chakraborty:Khare:2017}. Here we use it to illustrate convergences of the robit model  Markov chains.  Following \cite{Chakraborty:Khare:2017} we assign on the regression coefficient vector $\bb$ the Zellner's $g$-prior $\bb \sim \mathcal{N}_p(\bm{0}, \ g (X^T X)^{-1})$ with two choices of $g$: (a)  $g = 1000$ which induces a diffuse prior on $\bb$, and (b) $g = 3.49$ which ensures the trace-class property of the probit DA algorithm \citep{Chakraborty:Khare:2017}. For this data set we thus collectively consider 12 Markov chains from 3 models, 2 priors, and 2 Markov chain types (DA or sandwich). All 12 chains are initiated at the the maximum likelihood estimates $\beta_0 = -1.778$ (intercept), $\beta_1 = 4.374$, and $\beta_2 = 2.428$ obtained from the probit model. As the DA Markov chains for the Lupus data are known be slow-mixing \cite{Roy:Hobert:2007}, we  run each chain for  $10^6$ iterations, \textit{after} discarding the first $2 \times 10^6$ iterations as burn-in. We subsequently use the retained realizations from all the 12 Markov chains to compute (a) Markov chain autocorrelations up to lag 50, and (b) running means for each of the two non-intercept regression coefficients $\beta_1$ and $\beta_2$.  
       
% Although the main goal/contribution of this paper is theoretical, we provide illustrations to compare/contrast the performance of various relevant Markov chains on real data sets.

%\includegraphics[]{}

\begin{figure}[h]
    \centering
    \includegraphics[width = \linewidth]{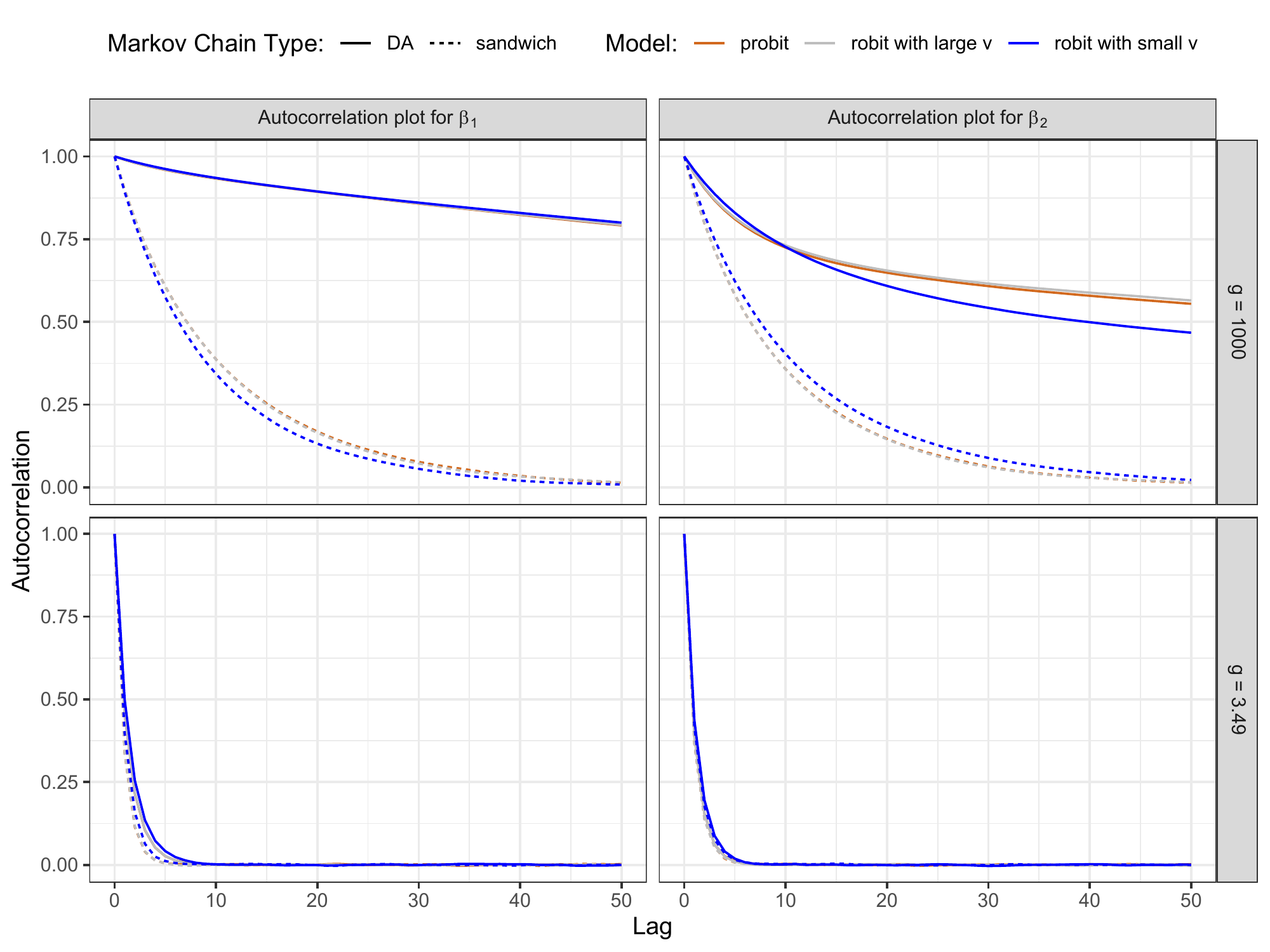}
    \caption{Autocorrelation plots for Markov chains run on the Lupus data set. }
    \label{fig:lupus_acf}
\end{figure}

\begin{figure}[h]
    \centering
    \includegraphics[width = \linewidth]{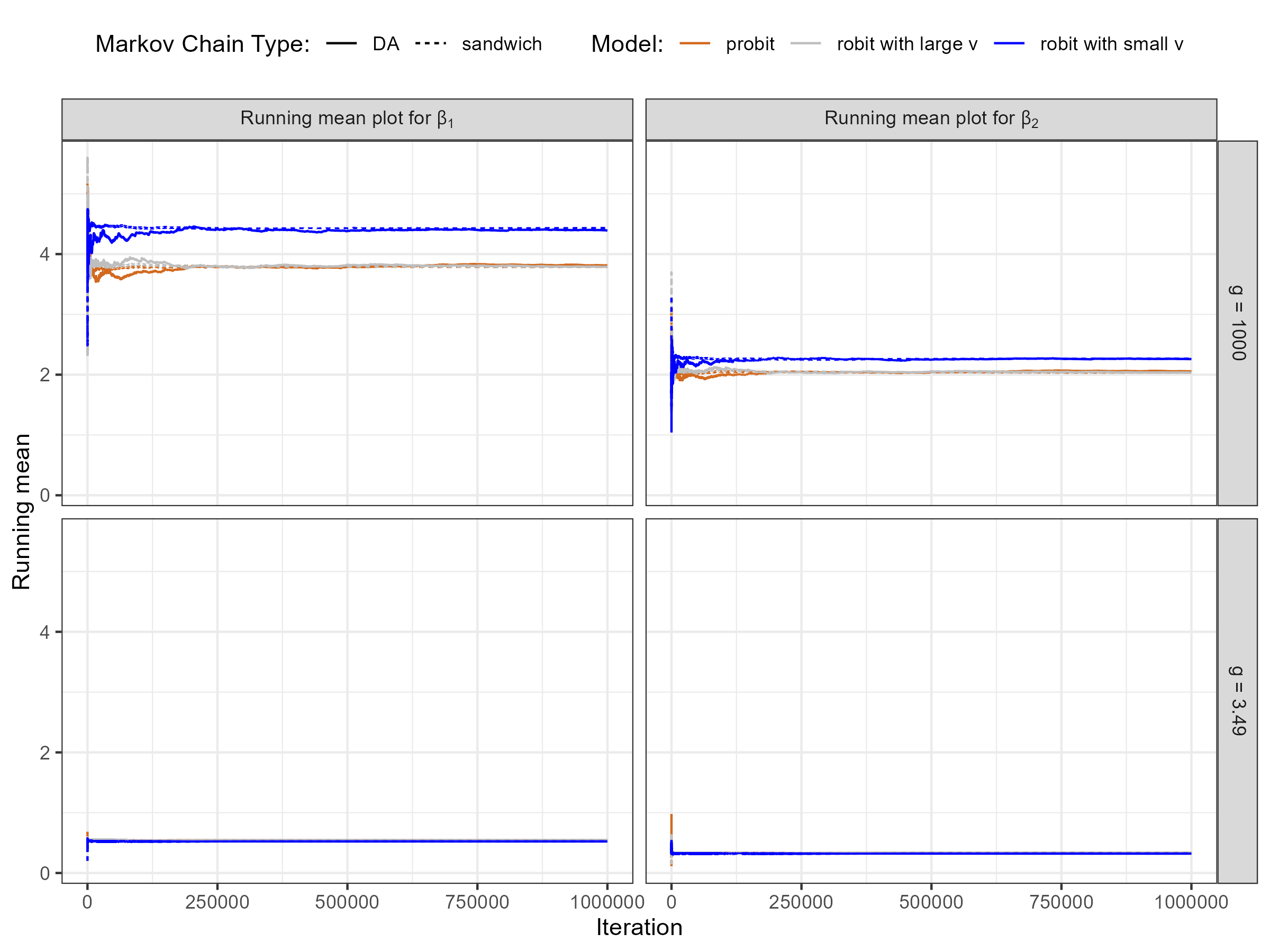}
    \caption{Running mean plots for the Markov chains run on the Lupus data set.}
    \label{fig:lupus_runmean}
\end{figure}

The autocorrelations and running means are displayed in Figures~\ref{fig:lupus_acf} and \ref{fig:lupus_runmean} as plot-matrices with the rows corresponding to priors ($g$ priors with $g = 1000$ and $g = 3.49$) and columns corresponding to $\bb$ components. Individual plots in the plot-matrices display as line diagrams autocorrelations ($y$-axis) plotted against lags ($x$-axis) in Figure~\ref{fig:lupus_acf}, or running means ($y$-axis) plotted against Markov chain iterations ($x$-axis) in Figure~\ref{fig:lupus_runmean}. A separate line is drawn for each model/Markov chain type combination (6 lines in total in each plot). The lines are color coded by models with the probit model, the robit model with small $\nu$, and the robit model with large $\nu$ being displayed as red, gray, and blue lines respectively. On the other hand, Markov chain types are displayed via line types: solid and dashed lines are used for DA and sandwich chains respectively. 

The following three observations are made from these two plots. First, for all three models sandwiching appears to aid substantial improvements in convergence and mixing over the original DA algorithm when the underlying prior is vague ($g = 1000$). This is demonstrated by both lowered autocorrelations in Figure~\ref{fig:lupus_acf} (top rows/panels) and stabler running means in Figure~\ref{fig:lupus_runmean} (top rows/panels) for the $\bb$-components in the sandwich chains. The improvements are not noticeable when the more informative prior with $g = 3.49$ is used. There, the DA and the sandwich chains display similar convergence properties, and the lines from the different model/Markov chain type combinations all effectively get superimposed at the displayed scale (bottom rows/panels in Figures~\ref{fig:lupus_acf} and \ref{fig:lupus_runmean}).  Interestingly, whereas the gain in sandwiching over the original DA algorithm is theoretically guaranteed for all choices of $g$ in the robit model (Theorem~\ref{thm1}), it has so far been theoretically guaranteed only for $g < 3.5$ in the probit model \cite{Chakraborty:Khare:2017}.
Second,  the running means from the DA chain are less stable than the sandwich chain; however,  they do converge to the same limit \textit{within the same model}. This is clearly expected since both the DA and sandwich chain have the same stationary distribution for $\bb$. This is particularly well-documented in the upper panels of Figure~\ref{fig:lupus_runmean}.  Third, when the vague prior ($g = 1000$) is used, the running means from the analogous chains for the  probit and the robit model with large $\nu$ become nearly identical with increasing iteration sizes, but they appear systematically different from the  running means for the robit model with small $\nu$. As noted before, this of course points to the relationship between a $t$-distribution with a large degrees of freedom and a standard normal distribution; as expected a robit model with a large $\nu$ gets well approximated by a probit model. When a more informative ($g = 3.45$) prior  is used the posteriors become less impacted by the systematic differences in the likelihoods and are more driven by the prior information; a fact well visualized in the lower panels of Figure~\ref{fig:lupus_runmean}. There, all chains from all models appear to share the same \textit{limiting} running means at the scale displayed.

\subsection{High Dimensional ($n < p$) Setting: Prostate Data Set}

\begin{figure}[h]
    \centering
    \includegraphics[width = \linewidth]{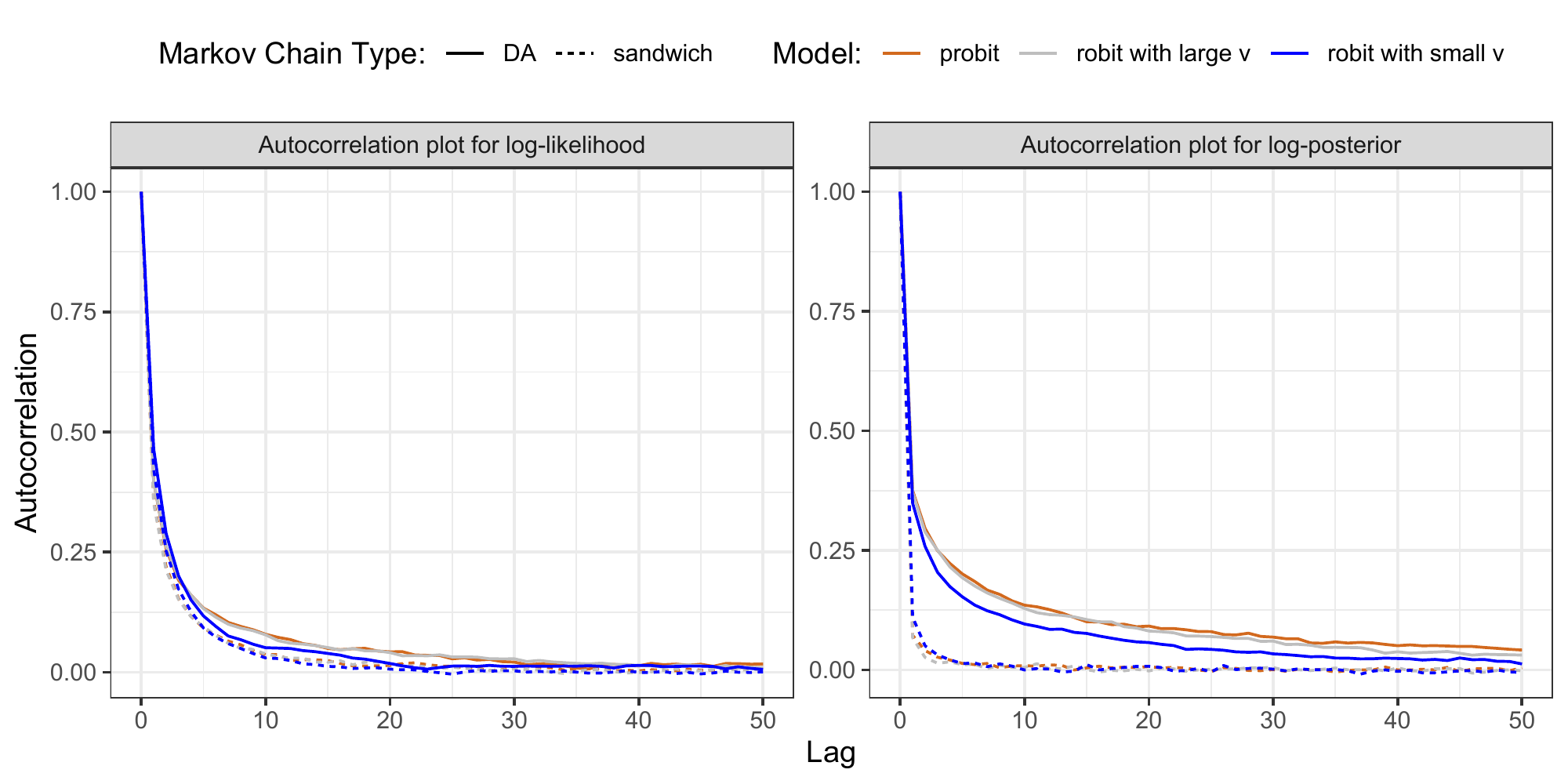}
    \caption{Autocorrelation plots for the Markov chains run on the Prostate data set.}
    \label{fig:prostate_acf}
\end{figure}

\begin{figure}[h]
    \centering
    \includegraphics[width = \linewidth]{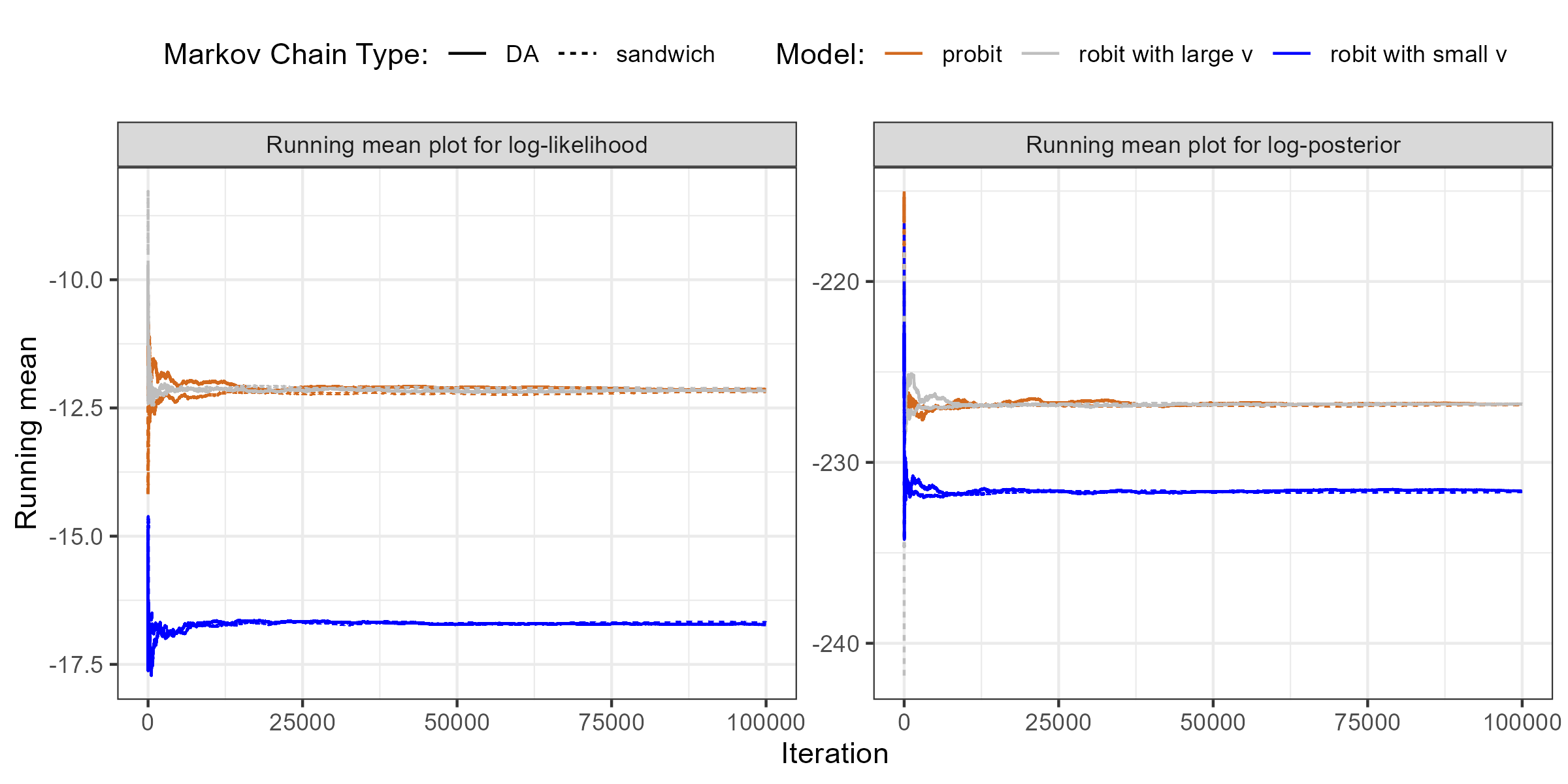}
    \caption{Running mean plots for the Markov chains run on the Prostate data set.}
    \label{fig:prostate_runmean}
\end{figure}

In the second example we consider the prostate cancer dataset from \cite{chunSparsePartialLeast2010}. The dataset records gene expressions of $50$ normal and $52$ prostate tumor samples at $6033$ arrays, of which we select the first $150$  arrays for our analysis. We are interested in regressing the binary cancer status (normal = 0, tumor = 1) on these selected 150 expression arrays (predictors). Similar to the analysis done in the previous section, we include an intercept term to the regression model to obtain the number of predictors $p = 151$ which is bigger than the total sample size of $n = 102$. We consider three models as before: (a) the probit model, (b) the robit model with a small $\nu \ (= 3)$, and (c) the robit model with a large $\nu \ (= 1000)$, and in each model assign \textit{independent standard normal priors} on the components of the regression coefficient vector $\bb$. For each model we then run two Markov chains -- the original DA chain, and the corresponding sandwich chain. All 6 chains are initiated at $\bb = \bm{0}$ and are run for $10^5$ iterations, \textit{after} discarding the first $2 \times 10^5$ iterations as burn-in. Subsequently, the (un-normalized) log-likelihood $\text{lik}(\bb)$ and (un-normalized) log-posterior density $\text{lpd}(\bb)$ values are calculated as univariate functions of $\bb$ on the retained realizations of each Markov chain. Here
\begin{align*}
\text{lik}(\bb) &= \sum_{i=1}^{n} \left\{ y_i \log F(\bm{x}_i^T \bb) + (1 - y_i) \log [1 - F(\bm{x}_i^T \bb)] \right\}, \ \text{and} \\ 
\text{lpd}(\bb) &= \text{lik}(\bb) - \frac{p}{2} \log(2\pi) - \frac12 \bb^T\bb,
\end{align*}
and $F$ is the normal/$t$ CDF associated with the probit/robit model. Finally for these computed log-likelihoods and log-posterior densities we calculate the Markov chain autocorrelations upto lag 50 and running means as done in the previous Lupus example. The resulting values are displayed in Figures~\ref{fig:prostate_acf} and \ref{fig:prostate_runmean} respectively. These figures follow the same color and line type conventions as considered in Figures~\ref{fig:lupus_acf} and \ref{fig:lupus_runmean}. 

The following observations are made from these figures. First, as expected, the sandwich chains (broken lines) are observed to have better convergences, i.e., smaller auto-correlations (Figure~\ref{fig:prostate_acf}) and stabler running means (Figure~\ref{fig:prostate_runmean}) than the original DA chains (solid lines). This is particularly well reflected among the log-posterior density values. Second, when focusing on the original DA chains, autocorrelations observed in both log-likelihood and log-posterior density values for the two robit chains are similar and are on average higher than the probit chain. By contrast all three sandwich chains appear to have similar auto-correlations. Third, for both the DA and sandwich chains, the ``limiting'' running means for the log-likelihood and log-posterior density values from the two robit models differ, but the probit model and the robit model with large $\nu$ appear similar. As noted before, this similarity in the likelihoods (and hence the posteriors) between the latter two models are attributable to the similarity between a $t$-distribution with large degrees of freedom and a standard normal distribution.

\appendix
%\section{Result related to matrix algebra}\label{appendix:matrix.algebra.results}
%\begin{result}
%	Let $A$ and $D$ be nonsingular matrices of orders $m$ and $n$ and $B$ be $m\times n$ matrix. Then
	%\begin{equation}
	%\left(A+BDB^T\right)^{-1} \ = \ A^{-1} \ - \ A^{-1}B\left(B^TA^{-1}B + D^{-1}\right)^{-1}B^TA^{-1}
	%\end{equation}
%\end{result}
%\noindent
%\begin{proof}
%See \cite{linear.stat.inference.book}.
%\end{proof}

\section{A Mill's ratio type result for Student's $t$ distribution}\label{appendix:mills.ratio.type.result}

\begin{lemma}\label{appendix:lemma.1}
    For $t>0$ we have
    \begin{equation}
        \frac{1}{\left(1 - F_{\nu}\left(t\right)\right)\left(t^2+\nu\right)^{\frac{\nu - 1}{2}}} \leq \frac{\sqrt{t^2+\nu - m}}{\kappa_m},
    \end{equation}
    where $\kappa_m = \Gamma\left(\left(\nu + 1\right)/2\right)\left(\nu - m\right)m^{\nu/2 - 1}/\left(2\sqrt{\pi}\Gamma\left(\nu/2\right)\right)$, $m$ is any arbitrary real number $\in \left(0, \nu\right)$ and as before $F_{\nu}(.)$ is the cdf of $t_{\nu}(0,1)$.
\end{lemma}
\vspace{0.02in}
\begin{proof}
Firstly, let us introduce the incomplete beta function ratio, defined by
\begin{equation}\label{eq:def.incomplete.beta.function.ratio}
    I_p\left(a, b\right) \ = \ \frac{1}{B(a,b)}\int_{0}^{p}u^{a-1}(1-u)^{b-1}du, \quad 0<p<1
\end{equation}
where, $a, b \in \R_+$, and $B(a,b) = \frac{\Gamma(a) \Gamma(b)}{\Gamma(a+b)}$.

Let, $f_{\nu}(.)$ be the pdf of $t_{\nu}(0,1)$. The cumulative distribution function $F_{\nu}(.)$ can be written in terms of $I$, the incomplete beta function ratio. From Chapter $28$ of \cite{continuous.univariate.distributions}, we have for $t>0$,
\begin{equation}
    F_{\nu}(t) \ = \ \int_{-\infty}^{t}f_{\nu}(u)du \ = \ 1 - \frac{1}{2}I_{x(t)}\left(\frac{\nu}{2}, \frac{1}{2}\right),
\end{equation}
where
\begin{equation}
    x(t) = \frac{\nu}{t^2+\nu}
\end{equation}
Now, from the definition \eqref{eq:def.incomplete.beta.function.ratio} of the incomplete beta function ratio, it follows that
\begin{align}
    I_{x(t)}\left(\frac{\nu}{2}, \frac{1}{2}\right) \ &= \ \frac{1}{B\left(\frac{\nu}{2}, \frac{1}{2}\right)} \ \int_{0}^{x(t)}u^{\frac{\nu}{2} - 1}\left(1-u\right)^{\frac{1}{2}-1}du\\
    %&\left[\text{where }B\left(x; a, b\right) \ \text{is the incomplete beta function defined as }\int_{0}^{x}u^{a-1}\left(1-u\right)^{b-1}du\right]\\
    &= \ \frac{\Gamma\left(\frac{\nu +1}{2}\right)}{\Gamma\left(\frac{\nu}{2}\right)\Gamma\left(\frac{1}{2}\right)} \ \int_{0}^{x(t)}u^{\frac{\nu}{2} - 1}\left(1-u\right)^{\frac{1}{2}-1}du\\
    &= \ \frac{\Gamma\left(\frac{\nu +1}{2}\right)}{\Gamma\left(\frac{\nu}{2}\right)\Gamma\left(\frac{1}{2}\right)} \ \int_{0}^{\frac{\nu}{t^2+\nu}}\frac{u^{\frac{\nu}{2} - 1}}{\sqrt{1-u}} \ du
    \label{eq:incomplete.beta.function.integral.TraceClass}
\end{align}
Again, if we consider the function $g(u) = u^{\nu/2 \ - 1}/\sqrt{1-u}$, then since
\begin{equation}
    g^{'}(u) = \frac{\left(\sqrt{1-u}\right)\left(\frac{\nu}{2}-1\right)u^{\frac{\nu}{2}-2} \ + \ u^{\frac{\nu}{2}-1}\frac{1}{2\sqrt{1-u}}}{1-u} > 0 \ \text{for all }u\in\left(0,1\right), \ \text{as }\nu>2
\end{equation}
which implies $g(u)$ is a strictly increasing function in $u$ for $u\in\left(0,1\right)$. Therefore, for the integral in the right hand side of \eqref{eq:incomplete.beta.function.integral.TraceClass}, we have
\begin{align}
    \int_{0}^{\frac{\nu}{t^2+\nu}}\frac{u^{\frac{\nu}{2} - 1}}{\sqrt{1-u}} \ du \ &= \ \int_{0}^{\frac{m}{t^2+\nu}}\frac{u^{\frac{\nu}{2} - 1}}{\sqrt{1-u}} \ du \ + \ \int_{\frac{m}{t^2+\nu}}^{\frac{\nu}{t^2+\nu}}\frac{u^{\frac{\nu}{2} - 1}}{\sqrt{1-u}} \ du\\
    & \ \left[\text{where }m \ \text{is any arbitrary real number in }\left(0,\nu\right)\right]\\
    &\geq \ \int_{\frac{m}{t^2+\nu}}^{\frac{\nu}{t^2+\nu}}\frac{u^{\frac{\nu}{2} - 1}}{\sqrt{1-u}} \ du\\
    &\geq \ \frac{\left(\frac{m}{t^2+\nu}\right)^{\frac{\nu}{2} - 1}}{\sqrt{1-\frac{m}{t^2+\nu}}}\int_{\frac{m}{t^2+\nu}}^{\frac{\nu}{t^2+\nu}} du\\
    &= \ \frac{\left(\frac{m}{t^2+\nu}\right)^{\frac{\nu}{2} - 1}}{\sqrt{1-\frac{m}{t^2+\nu}}} \ \times \ \frac{\left(\nu - m\right)}{t^2 + \nu}\\
    &= \ \frac{m^{\frac{\nu}{2}-1}\left(\nu - m\right)}{\sqrt{t^2 + \nu - m}} \ \times \ \frac{1}{\left(t^2+\nu\right)^{\frac{\nu - 1}{2}}}
    \label{eq:incomplete.beta.function.integral.LowerBound.TraceClass}
\end{align}
Hence, from \eqref{eq:incomplete.beta.function.integral.TraceClass} and \eqref{eq:incomplete.beta.function.integral.LowerBound.TraceClass}, we have
\begin{align}
    I_{x(t)}\left(\frac{\nu}{2}, \frac{1}{2}\right) \ &\geq \ \frac{\Gamma\left(\frac{\nu +1}{2}\right)m^{\nu/2 \ -1}\left(\nu - m\right)}{\sqrt{\pi}\Gamma\left(\frac{\nu}{2}\right)} \times \frac{1}{\sqrt{t^2+\nu - m}} \times \frac{1}{\left(t^2+\nu\right)^{\frac{\nu - 1}{2}}}\\
    \implies 1-F_{\nu}(t) \ &= \ \frac{1}{2}I_{x(t)}\left(\frac{\nu}{2}, \frac{1}{2}\right)\\
    &\geq \ \frac{\Gamma\left(\frac{\nu +1}{2}\right)m^{\nu/2 \ -1}\left(\nu - m\right)}{2\sqrt{\pi}\Gamma\left(\frac{\nu}{2}\right)} \times \frac{1}{\sqrt{t^2+\nu - m}} \times \frac{1}{\left(t^2+\nu\right)^{\frac{\nu - 1}{2}}}\\
    \implies \frac{1}{\left(1 - F_{\nu}\left(t\right)\right)\left(t^2+\nu\right)^{\frac{\nu - 1}{2}}}  \ &\leq \ \frac{\sqrt{t^2+\nu - m}}{\kappa_m}
\end{align}
where, $\kappa_m = \Gamma\left(\left(\nu + 1\right)/2\right)\left(\nu - m\right)m^{\nu/2 - 1}/\left(2\sqrt{\pi}\Gamma\left(\nu/2\right)\right)$, and $m$ is any arbitrary real number $\in \left(0, \nu\right)$.
\end{proof}

\begin{corollary}\label{appendix:corollary.1}
For $t>0$ we have
\begin{equation}
    \frac{1}{\left(1 - F_{\nu}\left(t\right)\right)} \leq \frac{\left(t^2+\nu\right)^{\frac{\nu}{2}}}{\kappa}
\end{equation}
where $\kappa = \Gamma\left(\left(\nu + 1\right)/2\right)\left(\nu - 1\right)/\left(2\sqrt{\pi}\Gamma\left(\nu/2\right)\right)$ and as before $F_{\nu}(.)$ is the cdf of $t_{\nu}(0,1)$.
\end{corollary}
\vspace{0.02in}
\begin{proof}
From Lemma \ref{appendix:lemma.1}, we know that for $t>0$ and for any arbitrary real number $m\in\left(0,\nu\right)$
\begin{equation}
    \frac{1}{\left(1 - F_{\nu}\left(t\right)\right)\left(t^2+\nu\right)^{\frac{\nu - 1}{2}}} \leq \frac{\sqrt{t^2+\nu - m}}{\kappa_m}
\end{equation}
Since $\nu > 2$, we can take $m=1$ to get the following 
\begin{align}
    \frac{1}{\left(1 - F_{\nu}\left(t\right)\right)\left(t^2+\nu\right)^{\frac{\nu - 1}{2}}} &\leq \frac{\sqrt{t^2+\nu - 1}}{\kappa}\\
    &\leq \frac{\sqrt{t^2+\nu}}{\kappa}\\
    \implies \frac{1}{\left(1 - F_{\nu}\left(t\right)\right)} &\leq \frac{\left(t^2+\nu\right)^{\frac{\nu}{2}}}{\kappa}
\end{align}
where, $\kappa = \Gamma\left(\left(\nu + 1\right)/2\right)\left(\nu - 1\right)/\left(2\sqrt{\pi}\Gamma\left(\nu/2\right)\right)$.
\end{proof}
%---------------------------------------------------------------

\end{document}